\documentclass[12pt]{amsart}
\usepackage{graphicx}
\usepackage{amssymb}
\usepackage{amsfonts}
\usepackage{amsmath}
\usepackage{mathtools}
\usepackage{array}
\usepackage{rotating}
\usepackage[matrix,frame]{xypic}

\usepackage{tikz-cd}
\usepackage{booktabs}

\usepackage{enumerate}
\usepackage[colorlinks,anchorcolor=black,citecolor=black,linkcolor=black]{hyperref}

\newcommand{\SP}{\mathcal{P}} 

\newcommand{\OSP}{\mathcal{OP}}
\newcommand{\Comp}{\mathcal{C}} 
\newcommand{\young}{\mathcal{Y}} 

\DeclareMathOperator{\IncAlg}{\mathfrak{I}}
\DeclareMathOperator{\linspan}{\mathcal{L}}
\DeclareMathOperator{\supp}{\mathrm{supp}}

\DeclareMathOperator{\Desc}{\mathcal{D}}

\newcommand{\SPleq}{\leq}
\newcommand{\SFleq}{\preceq}
\newcommand{\SFgeq}{\succeq}

\newcommand{\restr}{\rfloor}

\newlength{\tmplen}
\newcommand{\scriptpic}[1]{
  \setlength{\tmplen}{\unitlength}
  \setlength{\unitlength}{0.7pt}
  #1
  \setlength{\unitlength}{\tmplen}  
}

\usepackage{color}

\newcommand{\abs}[1]{
  \left\lvert
    #1
  \right\rvert}

\usepackage{svn}
\SVN $Revision: 91 $

\usepackage{marginnote}
\usepackage[a4paper,top=2cm,bottom=2cm,outer=3cm,marginparwidth=1cm]{geometry}

\theoremstyle{definition}
\newtheorem{example}{Example}[section]

\newtheorem{remark}[example]{Remark}
\newtheorem{definition}[example]{Definition}

\theoremstyle{theorem}
\newtheorem{theorem}[example]{Theorem}

\newtheorem{proposition}[example]{Proposition}

\newtheorem{lemma}[example]{Lemma}

\numberwithin{equation}{section}

\DeclareMathOperator{\PW}{\mathrm{PW}}

\def\Proof{\noindent \it Proof -- \rm}
\def\qed{\hspace{3.5mm} \hfill \vbox{\hrule height 3pt depth 2 pt width 2mm}
\bigskip}

\def\AA{{\mathcal A}}
\def\UU{{\mathcal U}}
\newcommand\alg[1]{\mathcal{#1}}
\newcommand\Exch{{\mathcal E}}
\newcommand\cS{{\mathcal S}}
\def\BB{{\mathcal B}}
\def\FQSym{{\bf FQSym}}
\def\WQSym{{\bf WQSym}}
\def\WSym{{\bf WSym}}

\newcommand\PiQSym{\mathit{\Pi\!QSym}}

\DeclareMathOperator{\ev}{\mathrm{ev}}
\DeclareMathOperator{\rgf}{\mathrm{rgf}}

\DeclareMathOperator\std{\mathrm{std}}

\def\Eul{{\mathcal E}}
\def\<{\langle}
\def\>{\rangle}

\def\IC{{\mathbb C}}

\def\N{{\bf N}}

\def\IN{\mathbb{N}}
\def\XX{{\mathbb X}}
\def\YY{{\mathbb Y}}

\def\K{{\bf K}}
\def\KK{\mathbb{K}}
\newcommand\IK{\KK}

\def\pack{\operatorname{pack}}

\def\F{{\bf F}}
\def\S{{\bf S}}

\def\M{{\bf M}}
\def\P{{\bf P}}

\def\E{{\mathbb E}}
\DeclareMathOperator{\IE}{{\mathbb E}}

\def\SG{{\mathfrak S}}
\def\H{{\mathcal H}}

\def\Sym{\mathbf{Sym}}
\def\sym{\mathit{Sym}}
\def\QSym{\mathit{QSym}}

\def\Des{\operatorname{Des}}

\def\m{{\bf m}}

\newcommand{\alf}[1]{\mathbb{#1}}
\newcommand{\alphA}{\alf{A}}
\newcommand{\alphX}{\alf{X}}
\newcommand{\alphY}{\alf{Y}}

\newcommand{\ii}{\underline{i}}
\newcommand{\oker}{\operatorname{\kappa}} 

\def\pev{\pi_{\rm ev}}

\def\gf#1#2{\left(\genfrac{}{}{0pt}{}{#1}{#2}\right)}
\def\ogf#1#2{\genfrac{}{}{0pt}{}{#1}{#2}}

\makeatletter
\newsavebox{\@brx}
\newcommand{\llangle}[1][]{\savebox{\@brx}{\(\m@th{#1\langle}\)}%
  \mathopen{\copy\@brx\kern-0.5\wd\@brx\usebox{\@brx}}}
\newcommand{\rrangle}[1][]{\savebox{\@brx}{\(\m@th{#1\rangle}\)}%
  \mathclose{\copy\@brx\kern-0.5\wd\@brx\usebox{\@brx}}}
\makeatother

\def\Tabvrule{\vrule width-0.4pt}       
\def\Tabhrule{\hrule \hrule height-0.4pt} 
\def\Tabstrut{\vrule height2.2ex 
                     depth0.8ex  
                     width0ex    
\relax}

\def\PasCase#1{\omit%
            $\vcenter{\hbox {\vbox to 0.4pt{}}
               \hbox{\makebox[3ex]{\Tabstrut$#1$}}}%
               \Tabvrule$}
\def\PasCasePoint{\PasCase{\cdot}}
\def\DessinCarre#1{%
    \vcenter{\hbox{}\hrule
             \hbox{\vrule\makebox[3ex]{\Tabstrut$#1$}\vrule}\Tabhrule}%
             \Tabvrule}
\def\GenRuban#1{\vcenter{\halign{&$\DessinCarre{##}$\cr#1}}\egroup}

\def\sTabvrule{\vrule width-0.4pt}
\def\sTabhrule{\hrule \hrule height-0.4pt}
\def\sTabstrut{\vrule height1.6ex depth0.6ex width0ex \relax}
\def\sDessinCarre#1{%
    \vcenter{\hbox{}\hrule
             \hbox{\vrule\makebox[2.3ex]%
                  {\sTabstrut$\scriptstyle#1$}\vrule}\sTabhrule}%
             \sTabvrule}
\def\sGenRuban#1{\vcenter{\halign{&$\sDessinCarre{##}$\cr#1}}\egroup}

\def\ruban{%
  \bgroup
  \let\ =\omit
  \let\\=\cr
  \let\x=\times
  \let\.=\PasCasePoint
  \offinterlineskip
  \GenRuban}

\def\sruban{%
  \bgroup
  \let\ =\omit
  \let\x=\times
  \let\\=\cr
  \offinterlineskip
  \sGenRuban}

\def\raff{{\rm Ref}}
\def\raffbis{{\rm WRef}}

\title{Combinatorial Hopf algebras in noncommutative probability}

\author[F. Lehner, J.-C.~Novelli, and J.-Y.~Thibon]%
{Franz Lehner, Jean-Christophe Novelli and Jean-Yves Thibon}

\address[Lehner]{
Institut f\"ur Diskrete Mathematik\\
TU Graz\\
Steyrergasse 30\\
A-8010 Graz, Austria}
\address[Novelli, Thibon] {Laboratoire d'Informatique  Gaspard Monge, Universit\'e Gustave Eiffel, CNRS, ENPC, ESIEE-Paris \\
5 Boulevard Descartes \\Champs-sur-Marne \\77454 Marne-la-Vall\'ee cedex 2 \\
France}
\email[Franz Lehner]{lehner@math.tugraz.at}
\email[Jean-Christophe Novelli]{novelli@univ-mlv.fr}
\email[Jean-Yves Thibon]{jyt@univ-mlv.fr} 

\keywords{Noncommutative symmetric functions, Eulerian polynomials, Eulerian idempotents}

\subjclass{05E05, 20C30, 60C05}

\date{Rev.\SVNRevision, \today}

\begin{document}

\begin{abstract}
  We extend the analogy between classical
  probability  and commutative symmetric functions
  to the noncommutative case. This involves the whole zoo of noncommutative
  symmetric functions   and we give a pedagogical overview the latter.
  We prove that the generalized moment-cumulant relations introduced  in
  \cite{HL} 
  are given by the action of the Eulerian idempotents on the Solomon-Tits algebras,
whose direct sum builds up the Hopf algebra of Word Quasi-Symmetric Functions $\WQSym$.
We prove $t$-analogues of these identities (in which the coefficient of $t$ gives back the original version),
and a similar $t$-analogue of Goldberg's formula for the coefficients of the Hausdorff series. This
amounts to the determination of the action of all the Eulerian idempotents on a
product of exponentials.

\end{abstract}

\maketitle

\tableofcontents{}


\section{Introduction}
The exponential function is a fundamental ingredient in many seemingly
unrelated areas of mathematics.
In the present paper we focus on two such areas
and establish a
deeply woven correspondence between noncommutative cumulants and
certain noncommutative generalizations of symmetric functions.
In both cases it is the noncommutativity which
leads to natural algebraic operations on set partitions and
give rise to interesting algebraic structures.
The starting point is the  observation
that the relation between  moments $m_n=\E(X^n)$ and classical cumulants
$\kappa_n(X)$ of a random variable $X$, encoded in the exponential generating functions
\begin{equation}
  \label{eq:defclass}
\sum_{n\ge0} m_n \frac{t^n}{n!} = \exp\left(\sum_{n\ge 1}\kappa_n \frac{t^n}{n!} \right)
\end{equation}
is, up to a rescaling by factorials, essentially the same as the relation
between the complete symmetric functions $h_n$ and the power sum symmetric functions
$p_n$
\begin{equation}
\label{eq:sumhntn=exprsumpn/ntn}
\sum_{n\ge 0}h_n t^n = \exp\left(\sum_{n\ge 1}\frac{p_n}{n} t^n\right).
\end{equation}
What seems to be a simple coincidence
actually goes much further:
cumulants are additive for independent random variables
\begin{equation}\label{eq:addcum}
\kappa_n(X+Y) = \kappa_n(X)+\kappa_n(Y)
\end{equation}
while the power sum symmetric functions are additive for unions of  alphabets
\begin{equation}
  \label{eq:addpn}
p_n(\alphX+\alphY) = p_n(\alphX) + p_n(\alphY).
\end{equation}
This analogy can be made rigorous in
terms of the Hopf algebra of symmetric functions
$\sym$. Indeed, the latter is freely generated (as an algebra)
both by the complete symmetric functions $h_n$ and the
power sums; and thus we can define a unique character by setting
\begin{equation}
\chi_X(h_n) = \frac{1}{n!}\E(X^n).
\end{equation}
Then the cumulants are given by
\begin{equation}
\kappa_n(X)=\chi_X\left(\frac{p_n}{(n-1)!}\right) 
\end{equation}
and the property that $\kappa_n(X+Y)=\kappa_n(X)+\kappa_n(Y)$ whenever $X$ and $Y$
are independent random variables corresponds in this context
to the fact that the power-sums are primitive elements for the coproduct
defined in \eqref{eq:Sym:coproduct}.
On a formal level this corresponds to Rota's umbral calculus
which was translated into the language of Hopf algebras early on 
\cite{Rota:foundationVIII,JoniRota:1979:coalgebras}.
The relations are summarized in Table~\ref{tab:class}.
\begin{table}
$$
    \begin{array}{r @{\hskip 1\tabcolsep} c @{\hskip 1\tabcolsep} lc|r @{\hskip 1\tabcolsep}c@{\hskip 1\tabcolsep}lc}
    \toprule{}
      \displaystyle
      \sum_{n\ge0} m_n \frac{t^n}{n!}
      &=&
      \displaystyle
            \exp\left(\sum_{n\ge 1}\kappa_n   \frac{t^n}{n!} \right)
  &
  \eqref{eq:defclass}
  &
      \displaystyle
    \sum_{n\ge 0}h_n t^n
      &=&
      \displaystyle
            \exp\left(\sum_{n\ge 1}\frac{p_n}{n} t^n\right)
      &
    \eqref{eq:sumhntn=exprsumpn/ntn}
  \\
      \midrule{}
      m_n & = &
      \displaystyle
      \sum_{\lambda\vdash n }
        \frac{n!}{z_\lambda } \kappa_\lambda
      &
        \eqref{eq:mn=sumkappapi}
      &
        h_n
      &=&
          \displaystyle
\sum_{\lambda\vdash n}\frac{p_\lambda}{z_\lambda}
  &
  \eqref{eq:hn=sumplambda/zlambda}
  \\
      \midrule{}
      \kappa_n(X+Y)
      &=&
            \kappa_n(X)+\kappa_n(Y)
      &
    \eqref{eq:addcum}
  &
    p_n(\alphX+\alphY)
      &=&
            p_n(\alphX) + p_n(\alphY)
   &
    \eqref{eq:addpn}
  \\
      \midrule{}
      \displaystyle
      m_n
      &=&
      \displaystyle
            \kappa_n + \sum_{k=1}^{n-1} \binom{n-1}{k-1}\kappa_km_{n-k}
      &
    \eqref{eq:classicalcumulantrecurrence}
  &
      \displaystyle
    nh_n
      &=&
      \displaystyle
            \sum_{k=1}^n p_k h_{n-k}
      &
        \eqref{eq:newtonrelation}
        \\
      \bottomrule  
    \end{array}
$$
  \caption{Correspondence between classical cumulants and symmetric functions}
  \label{tab:class}
\end{table}


It turns out that the multivariate setting leads to a better conceptual understanding.
Given a sequence $(X_i)_{i\ge 1}$ of random variables,
one can define multivariate moments
\begin{equation}
m_n := \E(X_1X_2\cdots X_n)
\end{equation}
as a multilinear functional.
Stochastic independence can then be algebraically reformulated
in terms of subalgebras of the algebra ${\mathcal X}$
of random variables.
A family $({\mathcal X}_i)_{i\in I}$ of subalgebras of  ${\mathcal X}$ 
is said to be \emph{independent} if for any choice of random variables
$X_1,X_2,\dots,X_n\in \bigcup \mathcal{X}_i$ the factorization
\begin{equation}\label{eq:mixmom}
\E(X_1X_2\cdots X_n) = \prod_{B\in\pi}\E\left(\prod_{i\in B}X_i\right)
\end{equation}
holds, where $\pi$ is any partition of $[n]$ such that for each block $B$
of $\pi$ the random variables $\{X_i \mid  i \in B\}$ come from the same
subalgebra ${\mathcal X}_{j(B)}$ and $j(B)\not=j(B')$ for $B\not=B'$.
The multivariate classical cumulants are then defined as
\begin{equation}
\kappa_n(X_1,X_2,\dots,X_n) = \left. \frac{\partial}{\partial t_1} \frac{\partial}{\partial
  t_2}\dotsm\frac{\partial}{\partial t_n}\right|_{{\bf t}=0} \IE e^{t_1X_1+t_2X_2+\dots+t_nX_n}
\end{equation}
which collapses to
  \eqref{eq:defclass} when $X_1=X_2=\dotsm=X_n=X$.
Multivariate
cumulants are multilinear maps $\kappa_n(X_1,X_2,\ldots,X_n)$ and exhibit the
fundamental property
of vanishing of mixed cumulants,
in the sense that $\kappa_n(X_1,X_2,\dots,X_n)=0$ whenever
at least two independent subalgebras occur.
The additive property \eqref{eq:addcum} is a simple consequence
of this.

Let us now turn to noncommutative probability where
the algebra of random variables is replaced by some noncommutative
algebra ${\mathcal A}$ and new notions of independence arise for which
the factorization of the moments \eqref{eq:mixmom} is replaced
by other identities. Each notion of independence 
gives rise to an appropriate version of cumulants, such that the vanishing of
mixed cumulants (or some weaker condition as in the case of monotone
independence) characterizes independence.

It is now a natural question whether the correspondence between classical
cumulants and symmetric functions can be lifted to
a correspondence between noncommutative cumulants and some noncommutative version of
symmetric functions. We will see that this is indeed the
case, but that it is actually the dual of the Hopf algebra of 
{\it quasi}-symmetric functions in noncommuting variables which hosts the correspondence.
This fact is concealed by the self-duality of the Hopf algebra of commutative
symmetric functions.

In an attempt to give a unified treatment of all these independences, the notions
of exchangeability system \cite{Leh04} and of spreadability system \cite{HL} have been introduced.
A spreadability system for a noncommutative probability space $({\mathcal A},\varphi)$
allows to produce independent copies $X^{(i)}$ of the random variables, and to formulate
independence in terms of identities satisfied by the moments 
$\varphi(X_1^{(i_1)}X_2^{(i_2)}\cdots X_n^{(i_n)})$.
It is in particular assumed that this quantity depends only on the \emph{packed word}
$u=\pack(i_1i_2\cdots i_n)$
(or equivalently \emph{ordered set partition} \cite{HL}) 
induced by the multiindex $i_1,i_2,\dots,i_n$, see §\ref{ssec:OSP} below.
This implies that for a spreadability  system, any choice of
a sequence $(X_i)$ of random variables determines a linear form on an appropriate Hopf
algebra based on packed words, namely, $\WQSym^*$, the graded dual of Word Quasi-Symmetric functions, 
also known as quasi-symmetric functions in noncommuting variables \cite{NT06,BZ09,NCSF7}.   
$\WQSym$ is a noncommutative version of the algebra of quasi-symmetric functions.
As shown by Hivert \cite{Hiv}, the quasi-symmetric polynomials are actually the invariants of
a very peculiar action of the symmetric group $\SG_n$ on polynomials in $n$ variables, called the quasi-symmetrizing
action, see §\ref{ssec:quasisym}.
This action can be extended to polynomials in noncommuting variables \cite{NCSF7}, and letting $n$ tend to infinity,
the algebra of invariants acquires a Hopf algebra structure, just as in the case of symmetric polynomials.

It turns out that the homogenous component of degree $n$ of its graded dual, $\WQSym_n^*$, can be identified
with the Solomon-Tits algebra of $\SG_n$ \cite {Sol76}, and that this identification is compatible with that
of $\Sym_n$, the space of noncommutative symmetric functions of degree $n$,
with the ordinary Solomon descent algebra.
That is, there is an embedding of Hopf algebras $\Sym\hookrightarrow \WQSym^*$ which is compatible with the internal products.
The different algebras are summarized in Table~\ref{tab:zoo}

\begin{table}
\begin{tikzcd}
  &
  \begin{matrix}
    \OSP{}\\
    13|7|246|5
  \end{matrix}
  \arrow[ld,"\ev" above] \arrow[rd,"\supp"] & \\
  \begin{matrix}
  \Comp{}\\
    (2,1,3,1)
  \end{matrix}
  \arrow[rd]& &
  \begin{matrix}
    \SP{}\\
    13/246/5/7
  \end{matrix}
  \arrow[ld]\\
  &
  \begin{matrix}
  \young{}\\
    (3,2,1,1)    
  \end{matrix}
\end{tikzcd}

\begin{tikzcd}
  &
  \begin{matrix}
    \WQSym{}\\
    \text{§\ref{ssec:wordquasisym}}
  \end{matrix}
 §\arrow[ld,"\theta" above,"\eqref{eq:WQSym:abelianization}" right]
 \arrow[leftarrow,rd,"\eqref{eq:WSyminWQSym}"]
 & \\
 \begin{matrix}
   \QSym{}
   \\
   \text{§\ref{ssec:quasisym}}
 \end{matrix}
 \arrow[leftarrow,rd,"\eqref{eq:SymintoQSym}"]
 & &
 \begin{matrix}
    \WSym{}\\
    \text{§\ref{ssec:WSym}}
 \end{matrix}
\arrow[ld,"\theta" above,"\eqref{eq:WSym:abelianization}" right]
 \\
 &
 \begin{matrix}
   \sym{} \\
   \text{§\ref{ssec:sym}}   
 \end{matrix}
\end{tikzcd}
\hfill{}
\begin{tikzcd}
  &
  \begin{matrix}
    \WQSym^*\\
    \text{§\ref{sssec:WQSym*}}
  \end{matrix}
 §\arrow[leftarrow,ld,"\theta^*" above,"\eqref{eq:embed:Symn->WQSym*}" right]
 \arrow[rd,"p" above, "\eqref{eq:p:WQsumtoWSym}" below]
 & \\
 \begin{matrix}
   \Sym{}
   \\
   \text{§\ref{ssec:Sym}}
 \end{matrix}
 \arrow[rd,"\eqref{eq:Sym->sym}"]
 & &
 \begin{matrix}
    \WSym^*\\
    \text{§\ref{sssec:WSym*}}
 \end{matrix}
\arrow[leftarrow,ld,"\theta^*" above,"\eqref{eq:theta*hlambda}" right]
 \\
 &
 \begin{matrix}
   \sym{} \\
   \text{§\ref{ssec:sym}}   
 \end{matrix}
\end{tikzcd}

  \caption{The zoo of noncommutative symmetric functions}
  \label{tab:zoo}
\end{table}


The moment-cumulant relations of \cite{HL} actually describe  the relation between the natural basis
$\N_u$ of $\WQSym^*$ and its internal product with the  Eulerian idempotents of the descent algebra $\Sym_n$.

\begin{table}
$$
    \begin{array}{c|cc|cc|c}
      \toprule{}
      & \multicolumn{2}{c|}{\text{monomials}}&\multicolumn{2}{c|}{\text{power
                                               sums}}&\text{index set}\\
    \midrule{}
      \sym & m_\lambda &   \eqref{eq:sym:mlambda}  & p_\lambda & \eqref{eq:def:pn} &\young\\
      \sym^* & h_\lambda &   \eqref{eq:def:hn}  & p_\lambda & \eqref{eq:def:pn} \\
      \midrule
      \QSym & M_I &   \eqref{eq:monomqsym} & &  &\Comp\\
      \Sym & S^I &   \eqref{Sn-G-fqs}  & \Phi_n & \eqref{def:Sym:Phi}\\      
      \midrule
      \WQSym &\M_u &  \eqref{eq:def:WQSym:Mu} & & &\OSP\\
      \WQSym^* & \N_u & \text{§\ref{sssec:WQSym*}} &\K_u &  \eqref{eq:def:WQSym*:Ku} & \\
      \midrule
      \WSym &  \m_\pi &   \eqref{eq:def:WSym:mpi} &   \Phi^\pi & \eqref{eq:WSym:Phi=summ}   &\SP    \\
      \WSym^* &N_\pi & \text{§\ref{sssec:WSym*}} & \phi_\pi &  \text{§\ref{sssec:WSym*}} \\
      \bottomrule
    \end{array}
$$
  \caption{A selection of bases of commutative and noncommutative symmetric functions and
    their duals}
\end{table}


\begin{table}
  $$
  \begin{array}{cc|c|cc|c}
    \toprule
    \multicolumn{2}{c|}{\text{lattice operation}}&\text{lattice}&\multicolumn{2}{c|}{\text{internal product}}&\text{algebra}\\
\midrule
    \pi\wedge\sigma &     \eqref{eq:piwedgerho} &\SP &N_\pi*N_\tau&    \eqref{eq:Npi*Ntau=Npiwedgetau}&\WSym^*
      \\
\midrule
    \pi \curlywedge \sigma &\eqref{eq:quasimeet}&\OSP &\N_u * \N_v & \eqref{eq:WQSym*:internalproduct} &\WQSym^*\\
\bottomrule
  \end{array}
  $$
  \caption{Lattice operations and internal products}
\end{table}


We provide simple conceptual proofs of most of the identities of \cite{HL}, and give a complete
description of the action of the Eulerian idempotents on $\WQSym^*$. As shown in \cite{FPT},
the celebrated Goldberg formula for the coefficients of the Campbell-Baker-Hausdorff series 
\begin{equation}
H(a_1,a_2,\ldots) = \log(e^{a_1}e^{a_2}\cdots) = \sum_{w\in A^*}g_w w,
\end{equation}
which is obtained by the action of the first Eulerian idempotent on the series
expansion of the product of exponentials \cite{Reu}, can be reduced to the calculation of
the action of the first Eulerian idempotent on the natural basis of the Solomon-Tits algebra $\WQSym^*$. 
Our results allow us to
compute the coefficients of the expansion
\begin{equation}
(e^{a_1}e^{a_2}\cdots)^t = e^{t\sum g_w w} = \sum_{w\in A^*}g_w(t) w 
\end{equation}
which reproduces the Hausdorff series in the coefficient of $t$.

Bialgebras and Hopf algebras have appeared in other approaches to
noncommutative probability as well, see for example
\cite{EbrahimiFardPatras:2016:splitting,FriedrichMcKay:2015:homogeneous,ManzelSchurmann:2016:noncommutative}
and the references therein.

\begin{table}[t]
  $$
    \begin{array}{r @{\hskip 1\tabcolsep} c @{\hskip 1\tabcolsep} lc|r @{\hskip 1\tabcolsep}c@{\hskip 1\tabcolsep}lc}  
      \toprule
      \multicolumn{4}{c|}{
    \text{$\Exch$-independence    
    }}
      &
      \multicolumn{4}{c}{        
      \text{internal product in $\WSym^*$ 
      }}
    \\
        \varphi_\pi(\underline{X})
    &=& \varphi_{\pi\wedge \ker(\underline{i})}(\underline{X})
    &\eqref{eq:E-indep}
&
   \hat\varphi_{\mathcal{X}}(N_\pi)
   &=&\hat\varphi_{\mathcal{X}}(N_\pi*N_\tau)
       &
    \eqref{eq:Exch:phiNpi=ph}                              
    \\
    %
    \varphi_\pi(\underline{X})
      &=&
    \displaystyle
      \sum_{\sigma\leq\pi}K_\sigma(\underline{X})
    &
    \eqref{eq:phipi=sumksig}
    &
      N_{\pi}
      &=&
    \displaystyle
          \sum_{\sigma\SFgeq\pi}\phi_{\sigma}
          &
      \eqref{eq:Nphi}
    \\
    %
        K_n(\underline{X}) 
      &=&
    \displaystyle
          \sum_{\pi\in\SP}
      \varphi_{\pi}(\underline{X})
  \,
  \mu(\pi,\hat{1}_n)
    &
    \eqref{eq:exch:Kn=sumphimu}
    &
            \phi_{[n] }
      &=&
    \displaystyle
          \sum_{\pi\in\SP}N_\pi\mu(\pi,\hat 1_n)
          &
      \eqref{eq:WSym*:phi=sumNmu}
    \\
    %
      K_n(\underline{X}^{(\underline{i})})
      &=& 0
    \text{ if } \ker\underline{i}\ne\hat{1}_n
          &
    \eqref{eq:Exch:mixedcumulantsvanish}
    &
      N_\pi* 
      \phi_{[n]}
      &=&0
      &
        \eqref{eq:prodNphi}
    \\
      \midrule
      \multicolumn{4}{c|}{              
    \text{
    $\cS$-independence 
    }}
    &
      \multicolumn{4}{c}{
      \text{
      internal product in $\WQSym^*$    
      }}
    \\
      \varphi_\pi(\underline{X})
      &=&
          \varphi_{\pi\curlywedge \oker(\underline{i})}(\underline{X})
          &
    \eqref{eq:S-indep}    
    &
      \hat\varphi_{\mathcal{X}}(\N_u)
      &=&
          \hat\varphi_{\mathcal{X}}(\N_u*\N_v)
          &
      \eqref{eq:Sindep}
      \\
      K_n(\underline{X})
      &=&
          \displaystyle
          \sum_{\sigma\in\OSP}
    \frac{(-1)^{\abs{\sigma}}}{\abs{\sigma}}
    \,
        \varphi_{\sigma}(\underline{X})
          &
    \eqref{eq:Kn:spread}
    &
      \K_u
      &=&
          \displaystyle
          \sum_{v\ge u}\frac{(-1)^{\ell(\ev(v))-\ell(\ev(u))}}{\ell(\ev(v),\ev(u))}\N_v
          &
      \eqref{eq:Ku=summuNv}
      \\
    \bottomrule
  \end{array}
  $$
  \caption{Independence and internal products}
\end{table}


\section{Words, formal power series and partitions}
The theory to be developed involves commutative and noncommutative formal power
series, integer partitions, set partitions,
ordered set partitions and various functions thereof.
In order to make the paper more accessible we will review these concepts in
considerable detail.

\subsection{Words, polynomials and formal power series}
An \emph{alphabet} is a finite or infinite set $\alphA$.
Its elements are called \emph{letters} or \emph{variables}.
An alphabet generates a free monoid $\alphA^*$ consisting of finite sequences
$w=a_1a_2\dotsm a_n$ of letters $a_i\in \alphA$, called \emph{words}
and the monoid multiplication is the concatenation of words.
Let $\IK$ be a field.
The multiplication extends linearly to the formal linear span
$\IK\langle \alphA\rangle=\linspan_\IK(\alphA^*)$ and turns it
into an algebra, the \emph{free associative algebra}
or algebra of \emph{polynomials in noncommuting variables}.
It is contained in the larger algebra of
\emph{noncommutative formal power series}
$\IK\llangle \alphA\rrangle$.
In the case of an infinite alphabet $\alphA$ we will still
denote by $\IK\langle \alphA\rangle$ the algebra of ``polynomials'',
i.e., in this case series with bounded degree.

If we let the variables commute, we obtain the usual algebras $\IK[\alphA]$
of (commutative) polynomials and $\IK[[\alphA]]$ of (commutative) formal power
series. 

\subsection{Spreadings and packed words}
\label{ssec:spreading}
Let $\alphA=\{a_1<a_2<\dotsm\}$ be an infinite totally ordered alphabet, usually we
will  assume that $\alphA=\IN$. 

A \emph{spreading} is an injective order preserving map $h:\alphA\to \alphA$
and the spreading of a word $w=w_1w_2\dotsm w_n$ is the multiplicative extension
$h(w)=h(w_1)h(w_2)\dots h(w_n)$ to the free monoid.

The inverse of a spreading is a \emph{packing}.
The \emph{packed word} $u=\pack(w)$ associated with a word $w\in \alphA^*$
is obtained by the following process:
If $b_1<b_2<\ldots <b_r$ are the different letters
occuring in $w$ rearranged in ascending order,
then $u$ is the image of $w$ by the homomorphism
$b_i\mapsto a_i$.
For example, $\pack(532725) = 321413$;
in other words, $\pack(w)$ is the lexicographically smallest word $u$ such that $w$ is a spreading of $u$.

A word $u$ is said to be \emph{packed} if $\pack(u)=u$. We denote by $\PW$ the
set of packed words.

\subsection{Integer partitions}
\label{ssec:integerpartitions}
A \emph{partition} of a positive integer $n$ is
a tuple of positive integers
$\lambda=(\lambda_1,\lambda_2,\dots,\lambda_k)$ such that
$\lambda_1\geq\lambda_2\geq\dots\geq\lambda_k$
and
\begin{equation}
  \label{eq:sigmalambda}
  \Sigma(\lambda) = \lambda_1+\lambda_2+\dotsm+\lambda_k = n
  .
\end{equation}
The shorthand notation is $\lambda\vdash n$.
\marginpar{$\vdash$}
The numbers $\lambda_i$ are called \emph{parts} of the partition
and the partial sums $s_i=\sum_{i=1}^j\lambda_j$ its \emph{descents}.
Sometimes it is more convenient to use the alternative notation
$$
\lambda = 1^{r_1}2^{r_2}\dots n^{r_n}
$$
where $r_i=r_i(\lambda)=\abs{\{j\mid \lambda_j=i\}}$ are the multiplicities.

The \emph{size} or \emph{length} of a partition $\lambda=(\lambda_1,\lambda_2,\dots,\lambda_k)$
is the number of its parts $\ell(\lambda) = k$.
\marginpar{ $\ell(\lambda)$}
We will denote the product of factorials of the parts
by
\marginpar{ $\lambda!$}
$$
\lambda! = \lambda_1!\lambda_2!\dotsm\lambda_k!
.
$$
For a sequence $(a_n)$ and a partition $\lambda=1^{r_1}2^{r_2}\dotsm$, we write
\begin{equation}
  \label{eq:alambda}
  a_\lambda = a_{\lambda_1}a_{\lambda_2}\dotsm a_{\lambda_k} = a_1^{r_1}a_2^{r_2}\dotsm
\end{equation}
Note however that this convention does \emph{not} apply to coefficients like
$\beta_\lambda$ and $z_\lambda$ defined below.

\subsection{Integer compositions}
\label{ssec:integercompositions}
A \emph{composition} of a positive integer $n$ is a tuple of positive integers
$I=(i_1,i_2,\dots,i_k)$ such that
\begin{equation}
  \Sigma(I) = i_1+i_2+\dotsm+i_k = n
  .
\end{equation}
The shorthand notation is $I\vDash n$.
\marginpar{$\vDash$}
Every composition has a unique nondecreasing rearrangement, which is a partition.
We usually denote compositions by capital roman letters $I,J,K,...$  to help
distinguish them from partitions,
but otherwise use the same terminology like length, descents, etc., as for partitions.
The set of integer compositions of $n$ is denoted by $\Comp_n$.

\subsection{Set partitions}
\label{ssec:setpartitions}
A partition of a finite set $X$ is a set of mutually disjoint subsets
$\pi=\{P_1,P_2,\dots,P_k\}$ such that $\bigcup P_i=X$.
The elements of a partition are called its \emph{blocks}.
Set partitions are in one-to-one correspondence with equivalence relations,
where we define $i\sim_\pi j$ if $i$ and $j$ are in the same block, i.e.,
the blocks are the equivalence classes.
Usually we deal with partitions of the set $X=[n]=\{1,2,\dots,n\}$,
which are called the set partitions of order $n$.
We denote the set of partitions of order $n$ by
\marginpar{$\SP_n$}
$\SP_n$.
The \emph{shape} of a set partition is the integer partition
$\lambda(\pi)$
\marginpar{$\lambda(\pi)$}
consisting of the sizes of the blocks of $\pi$.
We adapt the functions from §\ref{ssec:integerpartitions}
and write $\pi! = \lambda(\pi)!$, $\ell(\pi)=\ell(\lambda(\pi))$,
\begin{equation}
  \label{eq:a_pi=a_lambdapi}
a_\pi = a_{\lambda(\pi)}  
\end{equation}
etc.
We alert the reader however that convention  \eqref{eq:a_pi=a_lambdapi} will \emph{not}
apply to symmetric functions.
In the case of multiplicative functions
(see §\ref{sssec:reducedincalg})
summations over set partitions can be simplified to summations
over integer partitions weighted by the \emph{Bell coefficients}
which count 
the number of set partitions of given shape $\lambda$ 
\marginpar{$\beta_\lambda$}
\begin{equation}
  \label{eq:bellcoeff}
\beta_\lambda = 
\abs{\{\pi \in \SP_n\mid
  \lambda(\pi)=\lambda\}}=\frac{n!}{\lambda!r_1!r_2!\dotsm r_n!}
.
\end{equation}
\subsubsection{Kernel and rgf}
\label{sssec:kernel}
Set partitions arise naturally from the study of multiindices.
The \emph{kernel} of a function $f:X\to Y$,
denoted by $\ker f$,
\marginpar{$\ker$}
is the  partition of the set $X$ consisting of the equivalence classes of the equivalence
relation $i\sim j$ $\iff$ $f(i)=f(j)$
and similarly the kernel $\ker w$ of a multiindex or a word $w=a_1a_2\dotsm a_n$ arises
from the equivalence relation $i\sim j$ $\iff $ $a_i=a_j$.
Set partitions of order $n$ are in bijection with restricted growth functions (RGF)
of length $n$,
i.e., words $w_1w_2\dots w_n$ with letters $w_i\in\IN$ such that
$w_i\leq 1+\max\{w_1,w_2,\dots,w_{i-1}\}$.
Indeed every restricted growth function $w$ determines a partition
$\ker w$ and conversely,
if we order the blocks of a partition $\pi$
with respect to their smallest elements, the corresponding
restricted growth function
\marginpar{$\rgf$}
$\rgf(\pi)$ is obtained by assigning to $w_i$ the number of the block containing $i$.
For example, the partition $\pi = \{\{1,4,7\},\{2,8\},\{3\},\{5,6\}\}$,
in short notation $\pi=147/28/3/56$, is
mapped to the rgf $12314412$.

\subsubsection{Restriction, standardization, concatenation}
If $\pi=\{P_1,P_2,\dots,P_k\}$ is a partition of $X$
and $M\subseteq X$ is a subset,
we define the \emph{restriction}
$\pi\restr_M$
\marginpar{$\pi\restr_M$}
as the partition of $M$ consisting of the nonempty intersections $P_i\cap M$.

If $\pi$ is a partition of a finite totally ordered set $X=\{x_1<x_2\dots<x_n\}$
then the \emph{standardization} $\std\pi$\marginpar{$\std$}
is the partition $\pi'\in \P_n$
obtained by replacing the elements $x_i$ by their position $i$ in $X$.
For example, the partition $\pi = 147/28/3/56$ from
above restricted to the set $M=\{1,2,5,6,7\}$ results in
the partition
$\pi\restr_M = 17/2/56$
whose standardization is $\std(17/2/56) = 15/2/34$.

Two partitions $\pi\in\SP_m$ and $\rho\in\SP_n$
can be \emph{concatenated} to form an element
\marginpar{$\pi|\rho$}
$\pi|\rho\in\SP_{m+n}$.
To this end the blocks $\rho$ are added to $\pi$
after incrementing their entries by $m$ in oder to make them disjoint from $[m]$,
for example the concatenation of $\pi=15/2/34$ and $\rho=13/2/4$ is
$\pi|\rho=15/2/34/68/7/9$.

\subsubsection{Lattice structure}
The set partitions of fixed order $n$ form a lattice under
\emph{reverse refinement order}:
For partitions $\pi$ and $\sigma$ we write $\pi \SPleq \sigma$
$\sigma$ refines $\pi$, i.e.,
\marginpar{$\SPleq$}
\begin{equation}
  \label{eq:SPleq}
  \pi \SPleq \sigma
  \qquad\iff\qquad
  \forall P \in \pi\exists S \in \sigma:P \subseteq S
  .
\end{equation}

In other words, every block of $\sigma$ is a union of blocks of $\pi$
and the meet of two partitions $\pi$ and $\rho$
\begin{equation}
  \label{eq:piwedgerho}
\pi\wedge\rho = \{B\cap C \mid B\in\pi, C\in\rho\}\setminus\{\emptyset\}
\end{equation}
is the largest partition dominated by both $\pi$ and $\rho$.
The minimal element of this lattice is
$\hat{0}_n=\{\{1\},\{2\},\dots,\{n\}\}$ and the maximal element
$\hat{1}_n=\{\{1,2,\dots,n\}\}$.

Reverse refinement order is more natural from the point of view of
compatibility with the order on the Boolean lattice.
It is less so from the point of view of equivalence relations
and symmetric functions, where
traditionally \emph{refinement order} is preferred,
which we denote by $\SFleq$:
$\pi\SFleq\sigma$ means that $\pi$ is coarser than $\sigma$, i.e. that its
\marginpar{$\SFleq$}
blocks are unions of blocks of $\sigma$:
\begin{equation}
  \label{eq:SFleq}
  \pi \SFleq \sigma
  \qquad\iff\qquad
  \forall S \in \sigma \exists P \in \pi:S \subseteq P
  .
\end{equation}


\subsubsection{Incidence algebra and Möbius inversion}
The  incidence algebra of a partially ordered set $(P,\leq)$
was introduced by Rota \cite{Rota:foundationI}.
It is the algebra
$\IncAlg(P,\IK)$
\marginpar{$\IncAlg(P,\IK)$}
of 
functions supported on the set of pairs $\{(x,y)\in P\times P : x,y\in P; x\leq y\}$
with convolution product
$$
f*g(x,y) = \sum_{x\leq z\leq y} f(x,z)\, g(z,y)
.
$$ 
The algebra~$\IncAlg(P)$ has  the Kronecker function $\delta(x,y)$
as  unit element and a function $f\in\IncAlg(P)$
is invertible if and only if $f(x,x)$ is nonzero for every $x\in P$.
An example of an invertible function is the \emph{Zeta function},
which is defined as
\marginpar{$\zeta(x,y)$}
$\zeta(x,y)\equiv 1$.
Its inverse is called the \emph{M{\"o}bius function} of~$P$, 
denoted
\marginpar{$\mu(x,y)$}
$\mu(x,y)$.
For functions $F,G:P\to\IC$ we have the fundamental equivalence
 (``M{\"o}bius inversion formula'')
$$
  \forall x\in P:  F(x) = \sum_{y\leq x} G(y)
\qquad\iff\qquad
  \forall x\in P:  G(x) = \sum_{y\leq x} F(y)\,\mu(y,x)
$$
\subsubsection{Reduced incidence algebra and Möbius function}
\label{sssec:reducedincalg}
Before discussing the Möbius function in our setting we introduce
the reduced incidence algebra of the partition lattice,
which is one of the combinatorial incarnations of Fa\`a di Bruno's formula
\cite {Rota:foundationVI}.
First we note that any interval $[\pi,\sigma] = \{\rho \mid
\pi\SPleq\rho\SPleq\sigma \}\subseteq \SP_n$ is order isomorphic
to the direct product $\SP_1^{r_1}\times \SP_2^{r_2}\times\dotsm\times
\SP_p^{r_p}$,
where $p=\ell(\pi)$ and $\lambda(\sigma/\pi) = 1^{r_1}2^{r_2}\dots n^{r_n}\vdash p$
is the type of the partition $\sigma/\pi$ induced by $\sigma$ on the factor space $[n]/\pi$.
For example, let $\pi = 13/25/4/678$ and $\sigma = 134/25678$.
Then $\pi\leq\sigma$ and $\sigma/\pi = \{\{[13],[4]\},\{[25],[678]\}\}$.

Thus 
any sequence $(a_n)_{n \geq 1}$ gives rise to
\emph{multiplicative functions} on $\SP_n$, that  is an incidence function
$f\in\IncAlg(\SP_n)$ with values $f(\pi,\sigma) = a_{\lambda(\sigma/\pi)}$.
In particular, $f(0,\sigma) = a_{\lambda(\sigma)}$ and $f(\pi,1)=a_{\ell(\pi)}$
and 
the multiplicative functions form a subalgebra of $\IncAlg(\SP_n)$.
The convolution is given by
$
(a_n)*(b_n)=(c_n)
$
where
$$
c_n = \sum_{\pi\in\SP_n} a_\pi b_{\ell(\pi)}
= \sum_{\lambda\vdash n} \beta_\lambda a_\lambda b_{\ell(\lambda)}
$$
and $\beta_\lambda$ are the Bell coefficients  \eqref{eq:bellcoeff}.
In terms of generating functions this encodes
\emph{Fa\`a di Bruno's formula} 
(see \cite[Theorem~5.1.4]{StanleyVol2} or 
\cite[Proposition~5.9]{Aigner:1979})
which asserts that  the respective exponential generating functions
$F_a(z) = \sum_{k=1}^\infty \frac{a_k}{k!}z^k$,
$F_b(z) = \sum_{k=1}^\infty \frac{b_k}{k!}z^k$
and
$F_c(z) = \sum_{k=1}^\infty \frac{c_k}{k!}z^k$
satisfy the relation
\begin{equation}
  \label{eq:faadibrunocomposition}
  F_c(z) = F_b(F_a(z))
  .  
\end{equation}

The M\"obius function is multiplicative with characteristic sequence
$$
\mu_n = (-1)^{n-1}(n-1)!
$$
i.e, it evaluates to
\begin{equation}
  \label{eq:mu(pi,hat1)}
  \mu(\pi,\hat 1_n)=(-1)^{\ell(\pi)-1}(\ell(\pi)-1)!
\end{equation}
and we notice that the Bell coefficients
\eqref{eq:bellcoeff}
are related to  permutation statistics
by the identity
$$
\beta_\lambda 
=\frac{n!}{\abs{\mu_\lambda}z_\lambda}
.
$$
where the ubiquitous numbers
\marginpar{$z_\lambda$}
\begin{equation}
  \label{eq:zlambda}
  z_\lambda = r_1!r_2!\dotsm r_n!\,\lambda_1\lambda_2\dotsm\lambda_k,
\end{equation}
enumerate the centralizers of permutations
of cycle type $\lambda$ (cf.\ \cite[Prop.~2.2.2]{Meliot:2017:representation}).

\subsubsection{M\"obius algebra}
\label{sssec:moebiusalg}
The \emph{M\"obius algebra} \cite{Gre73}
of the partition lattice $\SP_n$
is the formal vector space spanned by partitions $\pi\in\SP_n$ with
the (commutative) multiplication
$\pi \cdot \rho=\pi\wedge \rho$.
It is semisimple and Solomon \cite{Sol67} showed that the elements
$$
\delta_\pi=\sum_{\rho\SPleq \pi} \rho
$$
form a complete set of orthogonal idempotents of the M\"obius algebra.
It is isomorphic to the dual of the Hopf algebra $\WQSym$ of symmetric
functions in noncommuting variables with the internal product, see §\ref{ssec:WSym}, in particular
  \eqref{eq:Npi*Ntau=Npiwedgetau}.

\subsection{Ordered set partitions}
\label{ssec:OSP}
An \emph{ordered set partition} or \emph{set composition} of order $n$
is a tuple  $(P_1,P_2,\dots,P_k)$ of nonempty disjoint subsets
$P_i\subseteq [n]$ such that $\bigcup P_i=[n]$;
that is, a set partition with a linear order on the blocks.
Similar to partitions we will denote ordered set partitions in shorthand notation
$P_1|P_2|\dots|P_k$, for example
$
(\{2,3\},\{4\},\{15\}) \simeq 23|4|15
$.
We denote the set of ordered set partitions of order $n$ by $\OSP_n$.
The ``support'' map
\marginpar{$\supp$}
\cite{Brown:2000:semigroups}
$\OSP_n\to\SP_n$ which ``forgets'' the order of blocks is
denoted by
\marginpar{$\bar{\pi}$}
$\supp:\pi\mapsto\bar{\pi}$,
for example
$\overline{23|4|15} = 15/23/4$.
\subsubsection{Kernels and packed words}
\label{sssec:orderedkernel}
All notations from set partitions naturally generalize to ordered set
partitions; the (ordered) \emph{kernel} $\oker(w)$ of a word $w=w_1w_2\dots w_n$
or function
\marginpar{$\oker(w)$}
$f:[n]\to T$ with values $f(i)=w_i$ in a totally ordered
set $T$
is obtained by ordering the blocks
of the unordered kernel $\ker f$
according to the values of the function $f$ on the individual
blocks. For example,
$\oker( 52235)=23|4|15$.
For any ordered set partition $\pi\in\OSP_n$ with $k$ blocks there is a unique word
$w=w_1w_2\dots w_n$ in letters $w_i\in[k]$ such that
$\pi=\oker(w)$. 
The word $w$ is called the \emph{packed word}
associated to $\pi$, denoted by $\pack(\pi)$, cf.\ §\ref{ssec:spreading}.
It follows that for any word $w$ we have $\pack(w) = \pack\ker(w)$.
In the example above, $\pack(23|4|15)=31123 = \pack(52235)$.

A packed word is a RGF if and only if the order of the blocks
of the corresponding partition is the natural one discussed in
§\ref{ssec:setpartitions}.

\subsubsection{Poset structure}
\label{ssec:OST:poset}
The reverse refinement order  \eqref{eq:SPleq} on set partitions is extended to
ordered set partitions as follows:
Two ordered set partitions
$\pi=(P_1, P_2,\dots,P_p), \sigma=(S_1, S_2,\dots,S_s) \in \OSP_n$,
satisfy  the order relation $\pi \leq \sigma$ if  the  following
conditions hold: 
\begin{enumerate}[(i)]
\item 
 $\bar{\pi}\leq\bar{\sigma}$ as set partitions. 
\item If $P_i \subseteq S_k, P_j \subseteq S_l$ for $i<j$, then $k \leq l$. 
\end{enumerate}

\subsubsection{Quasi-meet operation}
\label{sssec:quasimeet}
For lack of a unique minimal element
the ordered set partitions form a join-semilattice\footnote{Adjoining a minimal
  element turns it into a lattice isomorphic to the face lattice
  of the permutohedron
  \cite{BilleraSarangarajan:1996:combinatorics,Tonks:1997:relating}.
};
we will however replace the meet operation 
\eqref{eq:piwedgerho} with the \emph{quasi-meet operation}
of two ordered set partitions $\pi=(P_1,P_2,\dots,P_p)$ and
$\sigma=(S_1,S_2,\dots,S_s)$ given by
\marginpar{$\pi \curlywedge \sigma$}
\begin{equation}
  \label{eq:quasimeet}
\pi \curlywedge \sigma
=(P_1 \cap S_1, P_1\cap S_2 ,
  \ldots, P_1 \cap S_s, P_2 \cap S_1, P_2\cap S_2,\ldots, P_2\cap S_s, \dots,
  P_p \cap S_s)
\end{equation}
with empty intersections removed.
For example,
let  $\pi=1346|2|5$ and $\sigma=356|12|4$,
then $\pi\curlywedge\sigma=36|1|4|2|5$.

In terms of packed words
\begin{equation}
  \label{eq:ucurlywedgev}
 u\curlywedge v = \pack\binom{u}{v}
\end{equation}
where by $\binom{u}{v}$ we denote the word in {\it biletters} $\binom{u_i}{ v_i}$, lexicographically ordered
with priority to the top letter.
In the example above  $\pi=1346|2|5 = \oker(1\,2\,1\,1\,3\,1)$ and $\sigma=356|12|4 =
\oker(2\,2\,1\,3\,1\,1)$,
then 
the biletter word
$$
\binom{1\,2\,1\,1\,3\,1}{2\,2\,1\,3\,1\,1}
$$
consists of the biletters (lexicographically ordered)
$$
\binom{1}{1}
<
\binom{1}{2}
<
\binom{1}{3}
<
\binom{2}{2}
<
\binom{3}{1}
$$
and thus its packed word with respect to lexicographic order is
$$
\pack
\binom{1\,2\,1\,1\,3\,1}{2\,2\,1\,3\,1\,1}
=241351
$$
which corresponds to
$\pi\curlywedge\sigma=36|1|4|2|5$.

\subsubsection{Solomon-Tits algebra}
\label{sssec:solomontits}
The analog of the Möbius algebra from §\ref{sssec:moebiusalg}
is the Solomon-Tits algebra \cite{Sol76,PS2008}.
Indeed, the quasi-meet operation
  \eqref{eq:quasimeet} is associative and
  the Solomon-Tits algebra of order $n$ is the semigroup algebra
  spanned by the ordered set partitions.
  Moreover, the order forgetting map satisfies
  $\overline{\pi\curlywedge \sigma}=\bar{\pi}\wedge\bar{\sigma}$
  and therefore induces an algebra homomorphism
  from the Solomon-Tits algebra onto the Möbius algebra.
  In terms of packed words the quasi-meet operation corresponds to the internal product
  discussed in §\ref{ssec:wordquasisym} below.
  The Solomon-Tits algebra contains the Solomon descent algebra
  via the map
  \eqref{eq:embed:Symn->WQSym*} described below.

\section{Cumulants and independence}
\label{sec:cumulants}
\subsection{Classical cumulants}

The classical cumulants $\kappa_n$ of a random variable $X$ were
introduced by Thiele in 1889 via the exponential generating function
\begin{equation}
  \label{eq:defclasscumulants}
\sum_{n\ge 0}m_n\frac{t^n}{n!} = \exp\left(\sum_{n\ge 1}\kappa_n\frac{t^n}{n!}\right) 
\end{equation}
where $m_n=E[X^n]$ are the moments.
They are uniquely determined by the following properties
\begin{enumerate}[(K1)]
  
 \item \label{it:cum:K1}
  Additivity: If $X$ and $Y$ are independent random variables, then
    \begin{equation*}
      \kappa_n(X+Y) = \kappa_n(X)+\kappa_n(Y).
    \end{equation*}

   \item \label{it:cum:K2}
    Homogeneity: For any scalar~$\lambda$ the $n$-th cumulant is $n$-homogeneous:
    \begin{equation*}
      \label{eq:GeneralCumulants:axiom2}
      \kappa_n(\lambda X) = \lambda^n \kappa_n(X).
    \end{equation*}

   \item  \label{it:cum:K3}
    Universality:
     There exist universal polynomials $P_n$ in $n-1$ variables 
    without constant term such that
    \begin{equation*}
      \label{eq:GeneralCumulants:axiom3}
      m_n(X) = \kappa_n(X) + P_n(\kappa_1(X),\kappa_2(X),...,\kappa_{n-1}(X))
      .
    \end{equation*}
  \end{enumerate}

  Hasebe and Saigo  \cite{HS11b} showed 
  that property (K\ref{it:cum:K1}) can be replaced by the weaker condition
  \begin{enumerate}
   \item [(K1$^\prime$)] 
    Extensivity: If $X^{(i)}$, $i=1,2,\dots,N$, are i.i.d.\ copies of a random variables $X$, then
    \begin{equation*}
      \kappa_n(X^{(1)} + X^{(2)} + \dots+X^{(N)}) = N \kappa_n(X)
      .
    \end{equation*}
  \end{enumerate}

Combinatorially, cumulants and moments are related via
the identities
\begin{align}
  \label{eq:mn=sumkappapi}
  m_n &= \sum_{\pi\in\SP_n} \kappa_\pi = \sum_{\lambda\vdash n }
        \frac{n!}{z_\lambda } \kappa_\lambda
  \\
  \label{eq:kappan=summpi}
  \kappa_n&= \sum_{\pi\in\SP_n} m_\pi\mu(\pi,\hat{1}_n)
             = \sum_{\lambda\vdash n}
                (-1)^{\ell(\lambda)-1}({\ell(\lambda)-1})!\frac{n!}{z_\lambda }
            m_\lambda
\end{align}
and the recurrence
\cite{RS}
\begin{equation}
  \label{eq:classicalcumulantrecurrence}
  m_n = \kappa_n + \sum_{k=1}^{n-1} \binom{n-1}{k-1}\kappa_km_{n-k}
  .
\end{equation}

\subsection{Noncommutative probability and independence}
\label{ssec:independence}

A \emph{noncommutative probability space}, in short ncps, is a pair $(\AA,\varphi)$ where $\AA$ is a
unital algebra and $\varphi$ is a linear functional such that $\varphi(1)=1$.
The elements of $\AA$ are called  \emph{noncommutative random variables}.
More generally, one can replace the complex field by some unital algebra $\BB$ such that $\AA$ is a
$\BB$-module and consider
a unital linear 
map $\AA\rightarrow \BB$ to obtain an \emph{operator valued noncommutative  probability space}. 
Most constructions below carry over to the  operator valued setting, but this
is beyond the scope of the present paper.

Several notions of independence have been proposed,
among which free, Boolean and monotone independence
stand out, being the only ones satisfying certain natural axioms
\cite{Speicher:1997:universal,Muraki:2002:five,Muraki:2003:five}.
Each of these independences gives rise to a whole theory
and we will not review the specific details.
We rather concentrate on the common pattern
arising from an abstract invariance property
observed by the first named author in \cite{Leh04}.
More precisely, cumulants can be defined for any linear functional which is
invariant under the action of the symmetric group on its arguments,
and following a suggestion of D.~Voiculescu this setting was called
\emph{exchangeability system}.
Monotone independence is not covered in this setting,
as it is not invariant under permutations,
which have to be replaced by the weaker notion of spreadings.
To this end the theory was extended to the more general notion of \emph{spreadability systems}
in \cite{HL}.

\subsection{Exchangeability and spreadability}
\label{ssec:exchangeability} 

The main idea of \cite{Leh04,HL} is a reformulation of independence
in terms of the lattice operations   \eqref{eq:piwedgerho}
and   \eqref{eq:quasimeet}.
The starting point is  the observation that two classical random variables
$X$ and $Y$ are independent if and only if
the random vector $(X,Y)$ has the same distribution
as the random vector $(X^{(1)},Y^{(2)})$
where the vectors $(X^{(1)},Y^{(1)})$ and $(X^{(2)},Y^{(2)})$ are i.i.d.\ copies
of the vector $(X,Y)$. That is,
the joint distribution remains unchanged if 
the individual random variables are replaced with copies from two different
independent vectors.
The only ingredient needed for this is the invariance under
a certain action of the symmetric group
which will provide the connection to symmetric functions.

\subsubsection{Exchangeability}
\label{sssec:exchangeability} 

\begin{definition}
  Let  $(\AA,\varphi)$ be a noncommutative probability space.
  An \emph{exchangeability system} for $\AA$
  is a triple $\Exch=(\UU,\tilde\varphi,(I^{(i)})_{i\in \IN})$,
  where $(\UU,\tilde\varphi)$ is a noncommutative probablility space
  equipped with morphisms
  $I^{(i)}:\AA\rightarrow \UU$ such that $\varphi=\tilde\varphi\circ I^{(i)}$
  for all $i$, and, setting $X^{(j)}:=I^{(j)}(X)$, the variables
  $X^{(j)}$ are \emph{exchangeable}
  in the sense that
  \begin{equation}
    \tilde\varphi(X_1^{(i_1)}X_2^{(i_2)}\cdots X_n^{(i_n)})
    =
    \tilde\varphi(X_1^{(\sigma(i_1))}X_2^{(\sigma(i_2))}\cdots X_n^{(\sigma(i_n))})
  \end{equation}
  for any permutation $\sigma\in\SG_\infty$.
  The common value will be denoted by the multilinear functional
  $\varphi_\pi(X_1,X_2,\dots,X_n)$
  \marginpar{$\varphi_\pi$}
  where the kernel $\pi=\ker(\underline{i})\in\SP_n$ is the  set partition
  corresponding to the orbit of the index vector $(i_1,i_2,\dots,i_n)$
  under the symmetric group.
\end{definition}

For example, 
\begin{align*}
  \tilde\varphi(X_1^{(5)}X_2^{(3)} X_3^{(2)} X_4^{(7)} X_5^{(2)} X_6^{(5)})
  &= \tilde\varphi(X_1^{(1)}X_2^{(2)} X_3^{(4)} X_4^{(5)} X_5^{(4)}
    X_6^{(1)})\\
  &=:\varphi_{16/2/35/4}(X_1,X_2,X_3,X_4,X_5,X_6)
\end{align*}

\begin{definition}[$\Exch$-independence]
  \label{def:exchangeableindependence}   
   Let $\Exch=(\UU,\tilde{\varphi},(\iota^{(i)})_{i\ge1})$ be an exchangeability system for a ncps $(\alg{A},\varphi)$.  Let $I\subseteq \N$. 
\begin{enumerate}[(i)]
\item Subalgebras $(\alg{A}_i)_{i\in I}$ of $\alg{A}$ are said to be
  \emph{$\Exch$-independent} if for any tuple of indices $(i_1,i_2,\dots,i_n)\in I^n$, any tuple of random variables $(X_1,X_2,\dots,X_n)$ with $X_j\in \alg{A}_{i_j}$ and any set partition $\pi\in\SP_n$, we have  
  \begin{equation}
    \label{eq:E-indep}
  \varphi_\pi(X_1,X_2,\dots,X_n) = \varphi_{\pi\wedge \ker(i_1,i_2,\dots,i_n)}(X_1,X_2,\dots,X_n)
\end{equation}
where $\ker\underline{i}$ is the kernel set partition of the multindex
$\underline{i}$, cf.\ §\ref{sssec:kernel}.

\item Subsets  $(S_i)_{i \in I}$ of $\alg A$  are said to be
 \emph{$\Exch$-independent} if the algebras $\alg A_i$ generated by 
 $S_i$ are 
 $\Exch$-independent. In particular, a family
 $
 \bigl(
 (X_1(i),X_2(i),\dots,X_n(i))
 \bigr)_{i\in I}$
 of random vectors is said to be
 \emph{$\Exch$-independent} if the subalgebras $\alg A_i$
 generated by its respective entries  $\{X_1(i),X_2(i), \dots, X_n(i)\}$
 are $\Exch$-independent.  
\end{enumerate}
\end{definition}

That is, given a word in elements coming from independent
algebras, the sequence letter structure
For example, assuming $X$'s, $Y$'s and $Z$'s independent,
\begin{align*}
  \varphi_{1267/3458}(X_1Y_2X_3Z_4X_5Z_6X_7X_8)
  &=   \tilde{\varphi}(X_1^{(1)}Y_2^{(1)}X_3^{(2)}Z_4^{(2)}X_5^{(2)}Z_6^{(1)}X_7^{(1)}X_8^{(2)})
    \\
  &=   \tilde{\varphi}(X_1^{(1)}Y_2^{(2)}X_3^{(3)}Z_4^{(4)}X_5^{(3)}Z_6^{(5)}X_7^{(1)}X_8^{(3)})\\
  &=   \varphi_{1 7/2/358/4/6}(X_1Y_2X_3Z_4X_5Z_6X_7X_8)
\end{align*}

\subsubsection{Spreadability}
\label{sssec:spreadability} 

Spreadability generalizes exchangeability
in the same way that
quasisymmetric functions generalize symmetric functions,
to be discussed in the next section.
The main motivation for this generalization was Muraki's monotone independence
which does not satisfy exchangeability. 

\begin{definition}
  \label{def:spreadability}
A \emph{spreadability system} for $(\AA,\varphi)$ is a triple $(\UU,\tilde\varphi,(I^{(i)}))$,
where $(\UU,\tilde\varphi)$ is a $\BB$-valued noncommutative probablility space,
$I^{(i)}$ is a morphism $\AA\rightarrow \UU$ such that $\varphi=\tilde\varphi\circ I^{(i)}$
for all $i$, and, setting $X^{(j)}:=I^{(j)}(X)$, the variables
$X^{(j)}$ are \emph{spreadable} (or \emph{``quasi-exchangeable''})
in the sense that
\begin{equation}
\tilde\varphi(X_1^{(i_1)}X_2^{(i_2)}\cdots X_n^{(i_n)})
=
\tilde\varphi(X_1^{(h(i_1))}X_2^{(h(i_2))}\cdots X_n^{(h(i_n))})
\end{equation}
for any strictly increasing map $h:\IN\to\IN$;
equivalently, 
\begin{equation}
  \label{eq:phiii=phijj}
\tilde\varphi(X_1^{(i_1)}X_2^{(i_2)}\cdots X_n^{(i_n)})
=
\tilde\varphi(X_1^{(j_1)}X_2^{(j_2)}\cdots X_n^{(j_n)})
\end{equation}
whenever $\pack(i_1i_2\cdots i_n)=\pack(j_1j_2\cdots j_n)$.
Spreadably equivalent multiindices have common ordered kernel set
partition equal to  $\pi=\oker(\pack\underline{i})\in\OSP_n$ 
and the common value of the expectation
  \eqref{eq:phiii=phijj}
will be denoted by the multilinear functional
$\varphi_\pi(X_1,X_2,\dots,X_n)$.

More generally, we formally define 
 \begin{equation}
   \label{eq:phipixjij}
 \varphi_{\pi}(X_1^{(i_1)}, X_2^{(i_2)}, \ldots,X_n^{(i_n)}):=\varphi_{\pi \curlywedge \kappa(i_1,i_2,\ldots,i_n)}(X_1,X_2,\ldots,X_n),    
 \end{equation}
\end{definition}

For example, $\pack(532725) = 321413$ and spreadibility requires
\begin{align*}
  \tilde\varphi(X_1^{(5)}X_2^{(3)} X_3^{(2)} X_4^{(7)} X_5^{(2)} X_6^{(5)})
  &= \tilde\varphi(X_1^{(3)}X_2^{(2)} X_3^{(1)} X_4^{(4)} X_5^{(1)}
    X_6^{(3)})\\
  &=:\varphi_{35|2|16|4}(X_1,X_2,X_3,X_4,X_5,X_6)
\end{align*}

\begin{definition}[$\cS$-independence]
  \label{def:spreadableindependence}   
   Let $\cS=(\UU,\tilde{\varphi},(\iota^{(i)})_{i\ge1})$ be a spreadability
   system for a ncps $(\alg{A},\varphi)$.  Let $I\subseteq \N$. 
\begin{enumerate}[(i)]
\item Subalgebras $(\alg{A}_i)_{i\in I}$ of $\alg{A}$ are said to be
  \emph{$\cS$-independent} if for any tuple of indices $(i_1,i_2,\dots,i_n)\in I^n$, any tuple of random variables $(X_1,X_2,\dots,X_n)$ with $X_j\in \alg{A}_{i_j}$ and any set partition $\pi\in\SP_n$, we have  
  \begin{equation}
    \label{eq:S-indep}
  \varphi_\pi(X_1,X_2,\dots,X_n) = \varphi_{\pi\curlywedge \oker(i_1,i_2,\dots,i_n)}(X_1,X_2,\dots,X_n)
\end{equation}
where $\oker(\underline{i})$ is the ordered kernel set partition of the multindex
$\underline{i}$, cf.\ §\ref{sssec:orderedkernel}.

\item Subsets  $(S_i)_{i \in I}$ of $\alg A$  are said to be
 \emph{$\cS$-independent} if the algebras $\alg A_i$ generated by 
 $S_i$ are 
 $\cS$-independent. In particular, a family
 $
 \bigl(
 (X_1(i),X_2(i),\dots,X_n(i))
 \bigr)_{i\in I}$
 of random vectors is said to be
 \emph{$\cS$-independent} if the subalgebras $\alg A_i$
 generated by its respective entries  $\{X_1(i),X_2(i), \dots, X_n(i)\}$
 are $\cS$-independent.  
\end{enumerate}
\end{definition}
For example, assuming $X$'s, $Y$'s and $Z$'s independent,
\begin{align*}
  \varphi_{1267|3458}(X_1Y_2X_3Z_4X_5Z_6X_7X_8)
  &=
    \tilde{\varphi}(X_1^{(1)}Y_2^{(1)}X_3^{(2)}Z_4^{(2)}X_5^{(2)}Z_6^{(1)}X_7^{(1)}X_8^{(2)})
    \\
  &=
    \tilde{\varphi}(X_1^{(1)}Y_2^{(2)}X_3^{(4)}Z_4^{(5)}X_5^{(4)}Z_6^{(3)}X_7^{(1)}X_8^{(4)})
    \\
  &=   \varphi_{17|2|6|358|4}(X_1Y_2X_3Z_4X_5Z_6X_7X_8)
\end{align*}

\subsection{Cumulants}

It was shown in  \cite{HS11b} that cumulants satisfying
the axioms (K1$^\prime$)--(K3) listed above also characterize
free, Boolean and monotone independence.
This can be generalized to the framework of spreadability systems constructed
in \cite{HL}.
\begin{definition}
  Let $\cS=(\UU,\tilde{\varphi},(\iota^{(i)})_{i\ge1})$ be a spreadability
  system for a ncps $(\alg{A},\varphi)$.
  \marginpar{$N.X$}
  \emph{Rota's dot operator} is defined as
  \begin{equation}
    \label{eq:rotadot}
    N.X = X^{(1)}+X^{(2)}+\cdots+X^{(N)}
  \end{equation}
  Then, taking into account   \eqref{eq:phipixjij},
  one can show that for an ordered set partition $\pi\in\OSP_n$
  with $k$ blocks the formal expectation
  $\varphi(N.X_1\cdot N.X_2\cdots\cdot N.X_n)$ can be expanded into a polynomial in
  $N$ of degree $n$ without constant term:
  \begin{equation}
    \label{eq:phiNX=NKnX}
    \varphi((N.X_1)\cdot (N.X_2)\cdots\cdot (N.X_n)) = N K_n(X_1,X_2,\dots,X_n) + \alg{O}(N^2)
  \end{equation}
  The linear coefficient $ K_n(X_1,X_2,\dots,X_n)$ is a multilinear functional
  called the \emph{cumulant of order $n$}.
  More generally,
  $\varphi_{\pi}(N.X_1,N.X_2,\ldots,N.X_n)$ can be expanded into a polynomial in
  $N$ of degree $n$ divisible by $N^k$:
  \begin{equation*}
    \varphi_{\pi}((N.X_1),(N.X_2),\ldots,(N.X_n)) =
    N^k K_\pi(X_1,X_2,\dots,X_n)
    + \alg{O}(N^{k+1})
  \end{equation*}
  and the  lowest order coefficient $ K_\pi(X_1,X_2,\dots,X_n)$ is
  called \emph{partitioned cumulant}.
\end{definition}
\begin{example}
\begin{align*}
  \tilde{\varphi}((N.X)\cdot (N.Y))
  &= \sum_{i,j} \tilde{\varphi}(X^{(i)}Y^{(j)})\\
  &= \sum_{i} \tilde{\varphi}(X^{(i)}Y^{(i)})
    + \sum_{i< j} \tilde{\varphi}(X^{(i)}Y^{(i)})
        + \sum_{i>j} \tilde{\varphi}(X^{(i)}Y^{(i)})
  \\
  &= N\varphi_{11}(X,Y)
    \frac{N(N-1)}{2}(\varphi_{12}(X,Y)+ \varphi_{21}(X,Y)
  \\
  &=
    \begin{multlined}[t]
N(\varphi_{11}(X,Y) - \frac{1}{2}\varphi_{12}(X,Y) -
    \frac{1}{2}\varphi_{21}(X,Y))
    \\
    +    \frac{N^2}{2}( \frac{1}{2}\varphi_{12}(X,Y) + \frac{1}{2}\varphi_{21}(X,Y))
    \end{multlined}
\end{align*}
and thus
$$
K_2(X,Y) = \varphi_{11}(X,Y) - \frac{1}{2}\varphi_{12}(X,Y) - \frac{1}{2}\varphi_{21}(X,Y).
$$
The general formula involves the Möbius function   \eqref{eq:mu(pi,hat1)}.
In the spreadable case one obtains
\begin{align}
  K_n(X_1,X_2,\dots,X_n)
  &=   \sum_{\sigma\in\OSP_n}
    \frac{\mu(\sigma,\hat{1})}{\abs{\sigma}!}
    \,
    \varphi_{\sigma}(X_1,X_2,\ldots,X_n)
    \label{eq:Kn:spread}
  \\
  \notag
  &=   \sum_{\sigma\in\OSP_n}
    \frac{(-1)^{\abs{\sigma}}}{\abs{\sigma}}
    \,
    \varphi_{\sigma}(X_1,X_2,\ldots,X_n)
    ,
\end{align}
cf.~\cite[Theorem~3.11]{HL}.
\end{example}

$\Exch$-independence    \eqref{eq:E-indep} is then equivalent
to the \emph{vanishing of mixed cumulants},
i.e., 
if for any tuple of indices $(i_1,i_2,\dots,i_n)\in I^n$, any tuple of random variables $(X_1,X_2,\dots,X_n)$ with $X_j\in \alg{A}_{i_j}$ and any set partition $\pi\in\SP_n$, we have  
\begin{equation}
  \label{eq:Exch:mixedcumulantsvanish}
    K_\pi(X_1,X_2,\dots,X_n) = 0
    \text{ unless } \pi\leq \ker(i_1,i_2,\dots,i_n).
\end{equation}

\begin{remark}
  \label{rem:goodcartier}
  In the case of exchangeability systems
  formula \eqref{eq:Kn:spread}  
  simplifies to
  \begin{equation}
    \label{eq:exch:Kn=sumphimu}
  K_n(X_1,X_2,\dots,X_n)
  =   \sum_{\sigma\in\SP_n}
  \mu(\sigma,\hat{1})
  \,
  \varphi_{\sigma}(X_1,X_2,\ldots,X_n)
  \end{equation}
  and this can be rewritten   
  via in a  formula due to   Good \cite{Good:1975:new} and Cartier
  \cite{Percus:1975:correlation,Simon:1993:statistical}
  (see also Example~\ref{ex:goodcartier})
  as follows
  Let $\omega$ be a primitive root of unity of order $n$
  and put $X_i^{(\omega)} = \sum_{k=1}^n X_i^{(k)}$,
  then
  \begin{equation}
    \label{eq:exch:omegacom}
    K_n(X_1,X_2,\dots,X_n) =
    \frac{1}{n}\varphi(X_1^{(\omega)}X_2^{(\omega)}\dotsm X_n^{(\omega)}).
  \end{equation}
  This formula served as  definition of cumulants in \cite[Definition~2.1]{Leh04}.
\end{remark}


\section{Symmetric functions}
\label{sec:symfun}
In this section, we recall the basic definitions of the various Hopf algebras involved in the sequel:
ordinary ($\sym$) and noncommutative ($\Sym$) symmetric functions, word quasi-symmetric functions
($\WQSym$) and word symmetric functions ($\WSym$),
see the diagram in Table~\ref{tab:zoo}.
All algebras are taken over some field $\KK$ of characteristic 0.
%

Given a sequence $(X_i)$ of random variables, the moments 
$\varphi(X_1X_2\cdots X_n)$ determine a character of $\Sym$. We shall see that the
the structures introduced in \cite{Leh04,HL} allow to extend this character to a linear form
on $\WQSym$ (in the case of spreadability systems) or  $\WSym$ (in the case of
exchangeability systems).

\subsection {Ordinary symmetric functions}
\label{ssec:sym}
Given  an infinite set of commuting variables $\alphX=\{x_i|i\ge 1\}$
there is a natural action of the symmetric group $\SG_\infty$
which permutes the variables
$\sigma\cdot x_i=x_{\sigma(i)}$ and extends multiplicatively to monomials
\begin{equation}
  \label{eq:symact}
  \sigma\cdot x_1^{k_1}x_2^{k_2}\dotsm x_n^{k_n}
  =x_{\sigma(1)}^{k_1}x_{\sigma(2)}^{k_2}\dotsm x_{\sigma(n)}^{k_n}
  .
\end{equation}

This action lifts to the algebra of formal
power series $\KK[[\alphX]]$ and the invariant elements of the latter,
i.e., all power series whose coefficients satisfy the identity
$$
[x_1^{k_1}x_2^{k_2}\dotsm x_n^{k_n}]f(\alphX)
=[x_{\sigma(1)}^{k_1}x_{\sigma(2)}^{k_2}\dotsm x_{\sigma(n)}^{k_n}]f(\alphX)
$$
for any permutation $\sigma$ and any family of exponents $(k_1,k_2,\dots,k_n)$,
form the \emph{algebra of symmetric functions} $\sym=\sym(\alphX)$.

There are several natural bases and generating sets for this algebra:
The \emph{monomial symmetric functions} are given by the symmetrizations
of the standard monomials, indexed by partitions $\lambda\vdash n$.
Thus for a partition
$\lambda=\lambda_1\lambda_2\dots\lambda_r\vdash n$ the symmetrized ``monomial'' $m_\lambda$
is defined by
\marginpar{$m_\lambda$}
\begin{equation}
  \label{eq:sym:mlambda}
m_\lambda(\alphX) = \sum_{\substack{i_1<i_2<\dots<i_r\\
    \alpha\sim\lambda}}x_{i_1}^{\alpha_1}x_{i_2}^{\alpha_2}\dotsm x_{i_r}^{\alpha_r}
\end{equation}
where the sum runs over all ordered multiindices $\underline{i}$ and
compositions $\alpha\sim\lambda$, i.e., rearrangements of $\lambda$;
for example, symmetrizing the monomial $x_1^3x_2^2$ we obtain
$$
m_{3,2} = \sum_{i\ne j} x_i^3 x_j^2
.
$$
Clearly $(m_\lambda)_{\lambda\vdash n, n\in\IN}$ together with $1$
form a basis of $\sym$.

Other bases can be constructed from the following generating sets.

\begin{enumerate}[1.]
 \item 
  The \emph{complete homogeneous symmetric functions}
  \marginpar{$h_n$}
  \begin{equation}
    \label{eq:def:hn}
h_n(\alphX)=\sum_{i_1\le i_2\le \dots\le i_n}x_{i_1}x_{i_2} \dotsm x_{i_n}
\end{equation} 
that is, is the sum of all monomials of total degree $n$.
\item 
The 
\emph{elementary symmetric functions}
  \marginpar{$e_n$}
  \begin{equation}
        \label{eq:def:en}
  e_n(\alphX)=\sum_{i_1<i_2<\dots<i_n}x_{i_n}x_{i_{n-1}} \dotsm x_{i_1}
\end{equation} 
are the sum of all products of $n$ distinct variables.
\item 
 The \emph{power sum symmetric functions}
   \marginpar{$p_n$}
\begin{equation}
  \label{eq:def:pn}
  p_n=\sum_i x_i^n.
\end{equation}
\end{enumerate}
These
are invariant under permutations of the variables
and each of these sequences freely generates  $\sym(\alphX)$
as a commutative algebra:
\begin{equation}
  \sym = \KK[h_1,h_2,\ldots] = \KK[e_1,e_2,\ldots] = \KK[p_1,p_2,\ldots]
  .
\end{equation}
In order to obtain vector space bases we also have to consider their products
and using the convenient notation  \eqref{eq:alambda}
we find that all three families $\{h_\lambda \mid \lambda\vdash n,n\in \IN_0\}$,
 $\{e_\lambda \mid \lambda\vdash n,n\in \IN_0\}$ and $\{p_\lambda \mid \lambda\vdash n,n\in \IN_0\}$
form vector space bases of $\sym$.

\subsection{Partition lattice approach}
Alternatively, following Doubilet
\cite{Doubilet:1972:foundationsVII}, for any set partition $\pi$
we can define
\begin{equation}
  \label{eq:doubilet}
\begin{aligned}
  m_\pi &= \sum_{\ker \ii{}\, = \pi} x_{i_1}x_{i_2}\dotsm x_{i_n}\\
  p_\pi &= \sum_{\ker\ii{}\, \geq \pi} x_{i_1}x_{i_2}\dotsm x_{i_n}\\
  e_\pi &= \sum_{\ker\ii{}\, \wedge \pi =\hat{0}_n} x_{i_1}x_{i_2}\dotsm x_{i_n} 
\end{aligned}
\end{equation}


Notice however that convention   \eqref{eq:a_pi=a_lambdapi} does  \emph{not}
apply to symmetric functions, we rather have the following identities
\begin{align*}
  m_\pi &= r_1!r_2!\dotsm r_n!\, m_\lambda\\
  h_\pi &= \lambda!h_\lambda\\
  p_\pi &= p_{\lambda} \\
  e_\pi &=\lambda!\,e_{\lambda}
\end{align*}
where $\lambda$ is the shape of $\pi$.
We will see later
in §\ref{ssec:WSym}
that set partitions are more natural in the context of  noncommutative symmetric functions and 
these multiplicities will disappear.

\subsection{Generating functions}
The respective generating functions of the symmetric functions are
\begin{align}
  \label{eq:genfunc:ek}
  E(t) &= \sum_{k=0}^\infty e_kt^k= \prod (1+tx_i)
  \\
  \label{eq:genfunc:hk}
  H(t) &= \sum_{k=0}^\infty h_kt^k = \prod \frac{1}{1-tx_i}
  \\
  \label{eq:genfunc:pk}
  P(t) &= \sum_{k=1}^\infty \frac{p_k}{k}\,t^k = \sum\log\frac{1}{1-tx_i} 
\end{align}
and it follows that
\begin{equation}
  \label{eq:E-t*Ht=1}
  E(-t)H(t) = 1
\end{equation}
and that the power sums are related to the complete homogeneous
symmetric functions by the identity   \eqref{eq:sumhntn=exprsumpn/ntn}.

\subsection{Coproduct}
In addition to the algebra structure inherited from
the polynomial algebra,
there are coalgebra structures on $\sym$  which are
induced by set theoretic operations on the alphabet.
Given a second alphabet $\alphY$ disjoint from $\alphX$,
we  denote by $\alphX+\alphY$  the (still countable) union $\alphX\cup \alphY$
and we identify $\IC[[\alphX]]\otimes \IC[[\alphX]]$ with $\IC[[\alphX,\alphY]]$, i.e.,
we identify a tensor $f(\alphX)\otimes g(\alphX)$
with the element $f(\alphX)\,g(\alphY)$.
With this identification we have a coproduct 
$$
\label{eq:Deltaf=f(X+Y)}
\Delta(f)=f(\alphX+\alphY)
$$
which turns $\sym$ into a graded Hopf algebra.
On the level of generators it is easy to verify the following rules:
The complete homogeneous symmetric functions  form a \emph{group-like family} \cite[p.105]{CartierPatras:2021}
\begin{equation}
  \label{eq:Sym:coproduct}
  \Delta h_n = \sum_{k=0}^nh_k\otimes h_{n-k}
  .
\end{equation}
while the power sums
\eqref{eq:def:pn}
are \emph{primitive} elements, i.e.,
\begin{equation}
  \label{eq:Deltapn-primitive}
  \Delta p_n = p_n\otimes 1 + 1\otimes p_n
  .
\end{equation}


\subsection{Duality}
\label{ssec:duality}
The dual of a Hopf algebra $(H,\cdot,\Delta)$ has a natural Hopf algebra structure
$(H^*,*,\tilde{\Delta})$ such that
\begin{align*}
  \langle f*g,x  \rangle
  &= \langle f\otimes g, \Delta(x)\rangle \\
  \langle f,x\cdot y  \rangle
  &= \langle \tilde{\Delta}(f), x\otimes y\rangle \\
\end{align*}
It can be shown that $\sym$ is self-dual for the \emph{Hall scalar product}
\cite[Prop.~2.20]{Meliot:2017:representation} 
for which the complete homogeneous basis $h_\lambda$
and the monomial basis $m_\lambda$ are dual to each other:
\begin{equation}
  \label{eq:<hlambda,mmu>=deltalambdamu}
\langle h_\lambda,m_\mu\rangle = \delta_{\lambda,\mu}
\end{equation}
while the power sums are self-dual up to a rescaling by the numbers
\eqref{eq:zlambda}
\begin{equation}
  \label{eq:dual:plambdapmu}
\langle p_\lambda,p_\mu\rangle = z_\lambda \delta_{\lambda,\mu}  
.
\end{equation}


\subsection{Relation to classical cumulants}
The relation
  \eqref{eq:sumhntn=exprsumpn/ntn}
differs  from 
the relation 
between classical cumulants and moments
\eqref{eq:defclasscumulants}
by a simple rescaling.
Identifying $m_n$ with $n!\,h_n$, the cumulants $\kappa_n$ become identified with
$(n-1)!\,p_n$.
Then
the moment-cumulant formula
\eqref{eq:mn=sumkappapi} in terms of symmetric functions reads \cite[p.~62]{Meliot:2017:representation}
\begin{equation}
  \label{eq:hn=sumplambda/zlambda}
  h_n = \sum_{\lambda\vdash n}\frac{p_\lambda}{z_\lambda}
  ;
\end{equation}
and
the recurrence
\eqref{eq:classicalcumulantrecurrence}
corresponds to the Newton relation \cite[p.~63]{Meliot:2017:representation}
\begin{equation}
  \label{eq:newtonrelation}
nh_n = \sum_{k=1}^n p_k h_{n-k}
.
\end{equation}
In Doubilet's formulation this reads
\cite[p.~410]{Doubilet:1972:foundationsVII}
$$
h_\pi = \frac{1}{\lambda(\pi)!}\sum_{\sigma\SPleq\pi} \abs{\mu(0,\sigma)}p_\sigma
.
$$
This correspondence is best understood in terms of \emph{characters},
i.e., multiplicative linear functionals.
Any character of the algebra  $\sym$ is uniquely determined
by assigning arbitrary values to a free generating set.
Thus
a random variable $X$ determines a character $\chi_X$ of the algebra of symmetric
functions, by defining
\begin{equation}
  \label{eq:def:chiXhn}
\chi_X(h_n) = \frac{m_n}{n!}\ \Longleftrightarrow\ \chi_X(p_n)=\frac{\kappa_n}{(n-1)!}. 
\end{equation}
The values $\chi_X(e_n)$ then would produce cumulants for ``Steinhaus convolution''
\cite[p.~506]{Steinhaus:1985:selected}.

\subsection{Virtual alphabets}
A convenient symbolic notation for characters of $\sym$ is that of virtual alphabets:
given a sequence of free algebraic generators $g_n$ of $\sym$ (such as $h_n,p_n,e_n,\ldots$),
one denotes $\chi(g_n)$ by $g_n(\XX)$, where the symbol $\XX$ is called the virtual
alphabet associated with $\chi$. Then, if $\eta$ is another character whose virtual alphabet
is $\YY$, the convolution $\chi\star\eta := (\chi\otimes\eta)\circ\Delta$ is given
by 
\begin{equation}
\langle\chi\star\eta,f\rangle=f(\XX+\YY).
\end{equation}
On can thus associate with each random variable $X$ a virtual alphabet $\XX$,
such that $m_n(X)=n!h_n(\XX)$. If $X$ and $Y$ are independent, then 

\begin{equation}
\langle\chi_{X+Y}, h_n\rangle=\frac{m_n(X+Y)}{n!}=h_n(\XX+\YY)=
\langle\chi_X\otimes\chi_Y, \Delta h_n\rangle,
\end{equation}
that is, $\chi_{X+Y}$ is the convolution of  $\chi_X$ and $\chi_Y$.
The additivity of cumulants of sums of independent variables corresponds to  the fact
that the power-sums are primitive elements   \eqref{eq:Deltapn-primitive}.
This can be seen as another incarnation of the Hopf-algebraic
version of the \emph{umbral calculus} \cite{JoniRota:1979:coalgebras}.

%

\subsection{Another coproduct}
We also define the \emph{product of alphabets}
$\alphX\cdot \alphY := \alphX\times \alphY$ and the corresponding coproduct
\cite {Gessel:1984:multipartite}
\begin{equation}
  \label{eq:def:deltafXY}
  \delta f(\alphX,\alphY) = f(\alphX\cdot \alphY)
  .
\end{equation}

For this coproduct the power sums are group-like
$$
p_n(\alphX\cdot \alphY) = p_n(\alphX)\,p_n(\alphY)
$$
and from the expansion
\eqref{eq:hn=sumplambda/zlambda}
we obtain the expression
\begin{equation}
  \label{eq:deltahXY=sumpXpY/z}
  \delta h_n(\alphX,\alphY) = \sum_{\lambda\vdash n} \frac{p_n(\alphX)\,p_n(\alphY)}{z_\lambda}
  .
\end{equation}
Summing this identity and comparing with the generating function
\eqref{eq:genfunc:hk} for the alphabet $\alphX\cdot\alphY$ at $t=1$
we obtain the \emph{Cauchy product identity}
\cite[Theorem~2.18]{Meliot:2017:representation}
$$
\prod \frac{1}{1-x_iy_j}
= \sum_\lambda\frac{p_\lambda(\alphX)\,p_\lambda(\alphY)}{z_\lambda}
.
$$

\subsection{Internal product}
\label{ssec:class:internal}

By duality
(§\ref{ssec:duality})
the coproduct   \eqref{eq:def:deltafXY}
gives rise to an \emph{internal product}
$*$
on each homogeneous component $\sym_n$
and it follows from
  \eqref{eq:dual:plambdapmu}
that
the elements $\frac{p_\mu}{z_\mu}$ are orthogonal idempotents
for this product.

\begin{example}[The formula of 
  Good and Cartier]
  \label{ex:goodcartier}
  Pick a primitive root of unity $\omega$ of order $n$ and let
  $\alphY=\Omega_n=\{1,\omega,\ldots,\omega^{n-1}\}$ be the alphabet of $n$-th roots
  of unity, i.e., the roots of the polynomial $x^n-1$.
  Using the factorization $x^n-1=(x-1)(1+x+\dots+x^{n-1})$ we see that all the roots $\omega^k$ 
  for $1\leq k\leq n-1$  are roots of the second factor and therefore
  $p_k(\Omega_n)=0$ for $1\leq k\leq n-1$, whereas $p_n(\Omega_n)=n$.
  Thus we have from
  \eqref{eq:deltahXY=sumpXpY/z}
\begin{equation}
h_n(\Omega_n \alphX)=
\sum_{\lambda\vdash n}\frac{p_\lambda(\Omega_n) \, p_\lambda(\alphX)}{z_\lambda}
= p_n(\alphX).
\end{equation}
The evaluation of the character   \eqref{eq:def:chiXhn} on both sides
yields the classical version of identity \eqref{eq:omegacom} 
\begin{equation}\label{eq:omegacom}
  \E[(X^{(1)} + \omega X^{(2)} +\cdots +\omega^{(n-1)}X^{(n)})^n]
  =
  n\kappa_n
  ,
\end{equation}
known as \emph{Good's formula} in the mathematics literature  \cite{Good:1975:new}
and \emph{Cartier's formula} for Ursell functions in the physics literature
\cite{Percus:1975:correlation,Simon:1993:statistical}.
\end{example}

\subsection{Quasisymmetric functions}
\label{ssec:quasisym}
A weaker version of symmetry leads to the definition of quasisymmetric
functions \cite[§6.2]{Meliot:2017:representation}.
Given an infinite ordered alphabet $\alphX$,
a series $f(\alphX)\in \KK[\alphX]$ is called \emph{quasisymmetric} if for any increasing
families of indices $i_1<i_2<\dotsm<i_n$ and $j_1<j_2<\dotsm<j_n$ and any
family of exponents $k_1,k_2,\dots,k_n$ the
coefficients
of $x_{i_1}^{k_1}x_{i_2}^{k_2}\dotsm x_{i_n}^{k_n}$ and
$x_{j_1}^{k_1}x_{j_2}^{k_2}\dotsm x_{j_n}^{k_n}$ coincide.
Equivalently, a function $f(\alphX)\in \KK[\alphX]$ is quasisymmetric if it
is invariant under Hivert's  \emph{quasi-symmetrizing} action of the symmetric group
\cite[§3]{Hiv}, \cite{NCSF7} which acts on \emph{sets} of indices
\begin{equation}
  \label{eq:hivertaction}
  \sigma\cdot x_{i_1}^{k_1}x_{i_2}^{k_2}\dotsm x_{i_n}^{k_n}
  =x_{j_1}^{k_1}x_{j_2}^{k_2}\dotsm x_{i_n}^{k_n}
\end{equation}
where $\{j_1<j_2<\dotsm < j_n\}=\{\sigma(i_1),\sigma(i_2),\dots,\sigma(i_n)\}$.
For example, the transposition $\tau_{14}=(1\ 4)$ acts as follows:
\begin{align*}
  \tau_{14}\cdot x_1^3x_2^2x_3^4 &= x_2^3x_3^2x_4^4 \\
  \tau_{14}\cdot x_2^2x_3^4x_4^5 &= x_1^2x_2^4x_3^5 \\
  \tau_{14}\cdot x_1^3x_2^2x_3^4x_4^5 &= x_1^3x_2^2x_3^4x_4^5
\end{align*}
From these simple examples it is apparent
that contrary to the symmetrizing action   \eqref{eq:symact}
this action is not multiplicative,
i.i., $\sigma\cdot (ab) \ne (\sigma\cdot a) (\sigma\cdot b)$.
Yet the quasisymmetric functions form a subalgebra $\QSym\subseteq
\KK[\alphX]$, called the \emph{algebra of  quasisymmetric functions}.
A natural basis is given by the \emph{monomial quasisymmetric functions}
\marginpar{$M_I$}
\begin{equation}
  \label{eq:monomqsym}
  M_I
  = \sum_{j_1<j_2<\dotsm < j_r}
  x_{j_1}^{i_1}x_{j_2}^{i_2}\dotsm x_{j_r}^{i_r}
\end{equation}
indexed by compositions $I=i_1i_2\dotsm i_r$.
For example, restricting to the first 5 variables,
\begin{equation}
    \label{eq:M312}
  M_{3,1,2}=  \begin{multlined}[t]
   x_1^3x_2x_3^2
   + x_1^3x_2x_4^2
   + x_1^3x_3x_4^2
   + x_2^3x_3x_4^2
   + x_1^3x_2x_5^2
   \\
   + x_1^3x_3x_5^2
   + x_2^3x_3x_5^2
   + x_1^3x_4x_5^2
   + x_2^3x_4x_5^2
   + x_3^3x_4x_5^2
\end{multlined}
\end{equation}
Multiplication then gives rise to a quasi-shuffle product
\cite{Hoffman:quasi:2000}.

The algebra $\sym$ of symmetric functions is contained as a subalgebra via the
following 
embedding of the monomial symmetric functions
  \eqref{eq:sym:mlambda}
  \begin{equation}
    \label{eq:SymintoQSym}
  m_\lambda \mapsto\sum_{I\sim \lambda} M_I
\end{equation}
where the sum runs over all compositions which are permutations of the given
partition $\lambda$
\cite[p.~300]{Meliot:2017:representation};
for example,
$$
m_{2,2,1}\mapsto M_{2,2,1}+M_{2,1,2}+M_{1,2,2}
.
$$
The coproduct is again given by addition of alphabets,
this time however in a totally ordered manner
$\alphX+\alphY = (x_1,x_2,\dots,y_1,y_2,\dots)$ and consequently
it is obtained by deconcatenation of compositions
\begin{equation}
  \Delta M_I = \sum_{I'|I'' = I} M_{I'}\otimes M_{I''}
  ;
\end{equation}
for example,
$$
\Delta M_{2,1,3,2}
= 1\otimes M_{2,1,3,2}
  + M_2\otimes M_{1,3,2}
  + M_{2,1}\otimes M_{3,2}
  + M_{2,1,3}\otimes M_{2}
  + M_{2,1,3,2}\otimes 1
  .
$$

\section{Noncommutative symmetric functions}
\label{sec:ncsymfun}
The symmetric functions considered in §\ref{sec:symfun} have various
noncommutative generalizations.
Before turning to word symmetric functions
based on the natural analogues of Doubilet's formulas \eqref{eq:doubilet}
for noncommuting variables,
we start with a purely formal noncommutative generalization of the generating functions
\eqref{eq:genfunc:ek} and \eqref{eq:genfunc:hk} due to \cite{NCSF1},
see also \cite[p.~290]{Meliot:2017:representation}
and \cite{Thibon:2007:intro,Thibon:2011:intro} for a gentle introductions.

\subsection{The Hopf algebra $\Sym$ of noncommutative symmetric functions}
\label{ssec:Sym}
Let $\alphA=\{a_i|i\ge 1\}$ be an infinite totally ordered set of noncommuting variables.
We define the following analogues of     \eqref{eq:def:hn} 
\begin{equation}
\label{Sn-G-fqs}
S_n(\alphA) = \sum_{i_1\le i_2\le \dots\le i_n}a_{i_1}a_{i_2} \dotsm a_{i_n}
\end{equation} 
and      \eqref{eq:def:en} 
\begin{equation}
  \label{eq:def:Lambdan}
\Lambda_n(\alphA) = \sum_{i_1<i_2<\dots<i_n}a_{i_n}a_{i_{n-1}} \dotsm a_{i_1}
\end{equation} 
and call them respectively \emph{noncommutative complete symmetric functions} and
\emph{noncommutative elementary symmetric function}s of $\alphA$.
The reversed order in   \eqref{eq:def:Lambdan} is needed in order for 
the analogue of identity  \eqref{eq:E-t*Ht=1} to hold.

The $S_n(\alphA)$ freely generate a subalgebra of the noncommutative formal
power series over $\alphA$ which is denoted by
$\Sym(\alphA)$ and called \emph{noncommutative symmetric functions} \cite{NCSF1}.
It is isomorphic to the \emph{universal descent Hopf algebra}
\cite[§5.4] {CartierPatras:2021}
and should not be confused with the
algebra $\WSym$ of symmetric functions in noncommuting variables to be
considered below in §\ref{ssec:WSym}.

\subsubsection{Generating functions}
The $S_n(\alphA)$ and $\Lambda_n(\alphA)$ have the simple noncommutative generating
series
\begin{equation}
  \label{eq:def:sigmat}
\sigma_t(\alphA) := \sum_{n\ge 0}t^n S_n(\alphA)
             = \prod_{i\ge 1}^\rightarrow (1-ta_i)^{-1}
\end{equation}
and
\begin{equation}
\label{lambdaNC}
\lambda_{-t}(\alphA) := \sum_{n\ge 0}(-t)^n\Lambda_n(\alphA)
                 = \prod_{i\ge 1}^\leftarrow(1-ta_i)=\sigma_t(\alphA)^{-1}
\end{equation}
where we have set $S_0=\Lambda_0=1$, and $t$ is an indeterminate commuting
with $\alphA$;
making the variables commute yields the commutative generating series 
\eqref{eq:genfunc:ek}
and
\eqref{eq:genfunc:hk}.

We shall almost always forget about the alphabet $\alphA$ since the context is
generally clear.
As in the commutative case,
we can form a linear basis of the homogeneous component $\Sym_n$ from the
algebraic generators $S_n$ by assembling
\marginpar{$S^I$}
\begin{equation}
S^I = S_{i_1}S_{i_2}\cdots S_{i_r}
\end{equation}
for every 
integer composition $I$ of $n$.
The dimension of the homogeneous component $\Sym_n$ is thus $2^{n-1}$ for $n\ge 1$.
\subsubsection{Coproduct}

Let $\alf{A}$ and $\alf{B}$ be ordered alphabets and $\alf{A}\oplus\alf{B}$ their
ordinal sum.
Then clearly
\begin{align}
  \label{eq:Sym:sigma(A+B)}
  \sigma_t(\alf{A}\oplus\alf{B})
  &=\sigma_t(\alf{A})\sigma_t(\alf{B})
  &
    \lambda_t(\alf{A}\oplus\alf{B})
  &=\lambda_t(\alf{B})\lambda_t(\alf{A})
    .
\end{align}
Let now $\alf{A}+\alf{B}$ denote the union of the alphabets
such that the alphabets $\alf{A}$ and $\alf{B}$
mutually commute.
Then identifying the commuting copies with respectively $\alf{A}\otimes 1$ and
$1 \otimes \alf{A}$,
the coproduct is defined  similar to \eqref{eq:Deltaf=f(X+Y)} as
$$
\Delta f = f(\alf{A}+\alf{B})
$$
and endows $\Sym$ with the structure of a graded Hopf algebra
\cite[Proposition~5.5]{NCSF1}.
We infer from   
  \eqref{eq:Sym:sigma(A+B)}
that the generating function 
\eqref{eq:def:sigmat}
is group-like
$$
\Delta\sigma_t(\alf{A}) = \sigma_t(\alf{A})\otimes\sigma_t(\alf{A}) 
$$
and consequently the generators satisfy
\begin{equation}
  \label{eq:Sym:DeltaSn}
  \Delta S_n = \sum_{k=0}^n S_k\otimes S_{n-k}
\end{equation}
analogously to the commutative case \eqref{eq:Sym:coproduct}.
It follows from   \eqref{eq:Sym:DeltaSn} that the coproduct
is cocommutative and hence $\Sym$ cannot be self-dual. Its graded dual is
$QSym$, the algebra of quasi-symmetric functions
described above in
 §\ref{ssec:quasisym}. Similar to   \eqref{eq:<hlambda,mmu>=deltalambdamu}.
the dual basis of $S^I$ is the quasi-monomial basis
$M_I$ of $QSym$ from \eqref{eq:monomqsym},
see \cite[Theorem~6.17]{Meliot:2017:representation}.
The map
\begin{equation}
  \label{eq:Sym->sym}
  \begin{aligned}
    \Sym:\rightarrow \sym\\
    S_n:\mapsto h_n
  \end{aligned}
\end{equation}
extends to a morphism of Hopf algebras
dual to the natural embedding \eqref{eq:SymintoQSym} of  $\sym$ into $QSym$.

\subsubsection{Internal product}
\label{sssec:SYM:internalproduct}
Each homogeneous component $\Sym_n$ is endowed with an internal product $*$,
analogous to the classical one from §\ref{ssec:class:internal},
for which $\Sym_n$ is anti-isomophic to the Solomon descent algebra $\Sigma_n$
of $\SG_n$
(see \cite[Chapter~5]{NCSF1}, \cite[§6]{Thibon:2007:intro}, \cite {Schocker:2004}).
Recall that a \emph{descent} of a permutation $\sigma\in\SG_n$ is
a number $1\leq i\leq n-1$ such that $\sigma(i)>\sigma(i+1)$.
Denote this set by $\Desc(\sigma)$.
Solomon's descent algebra is the subalgebra of the group algebra
spanned by the elements
$D_S=\sum_{\Desc(\sigma)=S}\sigma$.
The anti-isomorphism mentioned identifies
the basis element $S^I$  with the formal sum
of all permutations whose descent composition is coarser than $I$.
In particular, the unit element of $\Sym_n$ is $S_n$, whose descent set is
empty, and thus formally
\begin{equation}
  \label{eq:sigma1*f=f}
  \sigma_1(\alphA)*f=f
\end{equation}
for any $f\in\Sym$.
This internal product is compatible with the coproduct in the sense that
the following ``splitting formula'' 
\cite[Proposition~5.2]{NCSF1}
\begin{equation}
  \label{eq:Sym:internal-coprod}
  (fg)*h = \mu_2( (f\otimes g) *_2 \Delta h).
\end{equation}
holds for arbitrary $f,g,h\in\Sym$.

\subsubsection{Solomon idempotents}
Finally, the primitive elements of any Hopf algebra span a Lie algebra
\cite[Lemma~3.1.1] {CartierPatras:2021}
which in the case of $\Sym$ is the free Lie algebra generated by
the  noncommutative power sum symmetric functions $\Phi_n$,
which are defined by analogy to \eqref{eq:sumhntn=exprsumpn/ntn} via the
generating function identity
\marginpar{$\Phi_n$}
\begin{equation}
  \label{def:Sym:Phi}
\sum_{n\ge 1}\frac{\Phi_n}{n} =\log\left(\sum_{n\ge 0}S_n\right).
\end{equation}
Then $\varphi_n:=\frac1n\Phi_n$ are idempotents for the internal product.
These are the \emph{Solomon idempotents}, also called \emph{(first) Eulerian idempotents}.

\subsubsection{Noncommutative cumulants as characters}
Given an infinite sequence $\mathcal{X}=(X_i)_{i\ge 1}$ of random variables
in the setting of §\ref{ssec:exchangeability},
the moment sequence
\begin{equation}
m_n := \varphi(X_1X_2\cdots X_n).
\end{equation}
determines a character $\hat\varphi_{\mathcal{X}}$ of $\Sym$ by setting
\begin{equation}
\hat\varphi_{\mathcal{X}}(S_n)=m_n
\end{equation}
(we drop the factorials which are irrelevant in this context).

The interpretation of the latter formalism
in terms of combinatorial Hopf algebras requires the introduction of $\WQSym$ (Word Quasi-Symmetric functions)
and its dual $\WQSym^*$. A spreadability system determines an extension of the character $\hat\varphi_{\mathcal{X}}$
to $\WQSym^*$, and the various notions of independence reflect certain symmetries of this extension.


\subsection{The algebra $\WQSym$ of word quasi-symmetric functions and its dual}

\label{ssec:wordquasisym}

The bases of the algebras considered in the following will be indexed
by packed words as defined in §\ref{ssec:spreading}.
We will usually represent the letters $a_i$ of the alphabet by $i$ in the
indexation of bases.

\subsubsection{Monomial basis}

To a packed word $u$ we associate the sum of all its spreadings, i.e., the
``polynomial''
\marginpar{$\M_u$}
\begin{equation}
  \label{eq:def:WQSym:Mu}
\M_u :=\sum_{\pack(w)=u}w\,.
\end{equation}
This is actually an infinite series,
but for example, restricting $\alphA$ to the first five integers, we obtain
\begin{multline*}
\M_{131132}= 131132 + 141142 + 141143 + 242243 \\
+ 151152 + 151153 + 252253 + 151154 + 252254 + 353354,
\end{multline*}
which, when ``written out'', represents the series
\begin{multline*}
  \M_{131132}=
  a_1a_3a_1^2a_3a_2 +
  a_1a_4a_1^2a_4a_2 +
  a_1a_4a_1^2a_4a_3
  +  a_2a_4a_2^2a_4a_3 
  +
  a_2a_5a_2^2a_5a_3
  \\
  +
  a_1a_5a_1^2a_5a_3 +
  a_2a_5a_2^2a_5a_3 +
  a_1a_5a_1^2a_5a_4 +
  a_2a_5a_2^2a_5a_4 +
  a_3a_5a_3^2a_5a_4
  .
\end{multline*}
For a word $w\in \alphA^*$ and a letter $a\in \alphA$ we denote by $\abs{w}_a$
\marginpar[ $\abs{\ }_a$]
the number of occurences
of the $a$ in $w$.
The \emph{evaluation} $\ev(w)$\marginpar{$\ev$}
is then the vector obtained
from the sequence $(\abs{w}_a)_{a\in \alphA}$ after removing all zeros,
which is a composition of $n$ where $n=\abs{w}$.
For example,
$$
\ev(131132) = (3,1,2)
.
$$
Equivalently, $ev(u)$ is the sequence of block sizes of the
corresponding ordered set partition $\kappa(u)$.

Under the abelianization
\begin{equation}
  \label{eq:WQSym:abelianization}
  \begin{aligned}
    \theta: \KK\langle \alphA\rangle& \rightarrow\KK[\alphX] \\
    M_u &\mapsto M_I
  \end{aligned}
\end{equation}
 the monomials $\M_u$ are mapped to the monomial    
quasi-symmetric function $M_I$ from 
 \eqref{eq:monomqsym}
where $I=(i_1,\ldots,i_r)=\ev(u)$ is the evaluation vector of $u$.
Thus
$$
\theta(\M_{131132}) = M_{3,1,2}
$$
which is shown in     \eqref{eq:M312}.
The $\M_u$ span a subalgebra 
of $\K\langle \alphA\rangle$, called $\WQSym$ for Word
Quasi-Symmetric functions,
consisting of the invariants of the noncommutative
version of Hivert's quasi-symmetrizing action
\eqref{eq:hivertaction}
which is defined by
$\sigma\cdot w = w'$ where $w'$ is such that $\std(w')=\std(w)$ and
$\theta(w')=\sigma\cdot\theta(w)$.
Hence, two words are in the same
$\SG(\alphA)$-orbit iff they have the same packed word.
For example, the transposition $\tau_{14}=(1\ 4)$ acts as follows:
\begin{align*}
  \tau_{14}\cdot 131132 &= 242243\\
  \tau_{14}\cdot 34324 &= 23213\\
  \tau_{14}\cdot 13141324 &= 13141324
  .
\end{align*}


\subsubsection{The internal product on the dual algebra $\WQSym^*$}
\label{sssec:WQSym*}
Let $\N_u\in\WQSym^*$
\marginpar{$\N_u$}
be the dual basis of $\M_u$, i.e.,
$\langle N_u,M_v\rangle = \delta_{u,v}$.
Define an internal product
on $\WQSym^*_n$ by \cite{NT06}
\begin{equation}
  \label{eq:WQSym*:internalproduct}
\N_u * \N_v = \N_{\pack{u\choose v}},
\end{equation}
where ${u\choose v}$ denotes the word in {\it biletters} ${u_i\choose v_i}$, lexicographically ordered
with priority to the top letter, cf.~§\ref{sssec:quasimeet}.
Then, 
\begin{proposition}
$(\WQSym^*,*)$ is anti-isomorphic to the Solomon-Tits algebra from §\ref{sssec:solomontits}.
\end{proposition}
Indeed, if $w=\pack{u\choose v}$
then $\oker(w) = \oker( u) \curlywedge\oker( v)$ defined in \eqref{eq:quasimeet}.


Finally the map
\begin{equation}
  \label{eq:embed:Symn->WQSym*}
  \begin{aligned}
\theta^*:\Sym_n&\to \WQSym^*\\
S^I& \mapsto \sum_{\ev(u)=I}\N_u
  \end{aligned}
\end{equation}
defines an embedding of Hopf algebras,
for example,
$$
S^{1,2}
\mapsto
\N_{122}+\N_{212}+\N_{221}
\qquad
S^{21}\mapsto\N_{112}+\N_{121}+\N_{211}
.
$$
This map is compatible with the internal product $*$, that is, it induces
the standard embedding of the descent algebra into the Solomon-Tits algebra.

\subsubsection{Encoding of spreadability systems and exchangeability systems}
From the definition of spreadability in §\ref{sssec:spreadability} 
we infer that every sequence  $\mathcal{X}=(X_i)_{i\ge 1}$ of random variables
in $\AA$
induces a linear map $\hat\varphi_{\mathcal{X}}:\ \WQSym^*\rightarrow \BB$
\marginpar{$\hat\varphi_{\mathcal{X}}$}
by
\begin{equation}
  \hat\varphi_{\mathcal{X}}(\N_u) := \varphi(X_1^{(u_1)}X_2^{(u_2)}\cdots X_n^{(u_n)}).
\end{equation}
In the notation of Definition~\ref{def:spreadability}
this is $\varphi_\pi(X_1X_2\cdots X_n)$, where $\pi=\oker( u)$
is the ordered set partition encoded by the packed word $u$, see §\ref{ssec:OSP}.

The notion of ${\mathcal S}$-independence
from Definition~\ref{def:spreadableindependence}   
 can be reformulated in terms of the internal product
 of $\WQSym^*$ as follows:
 A family of subalgebras $\AA_i$ is ${\mathcal S}$-independent
if
\begin{equation}
  \label{eq:Sindep}
\hat\varphi_{\mathcal{X}}(\N_u)=\hat\varphi_{\mathcal{X}}(\N_u*\N_v)
\end{equation}
for every word $u$ and every sequence $\mathcal{X}$ such that $X_j\in \AA_{v_j}$ for every $j$.

An exchangeability system (Definition~\ref{sssec:exchangeability})  is a spreadability system satisfying
\begin{equation}
\hat\varphi_{\mathcal{X}}(\N_{\sigma(u)}) = \hat\varphi_{\mathcal{X}}(\N_u)
\end{equation}
for every sequence $\mathcal{X}$ and 
for all permutations $\sigma$ of the alphabet of $u$. In this case,
$\hat\varphi_{\mathcal{X}}$ can be interpreted as a character of $\WSym^*$,
an algebra based on set partitions, to be defined below.

\subsection{The algebra $\WSym$ of symmetric functions in noncommuting variables}
\label{ssec:WSym}

Let $\alphA$ be
an alphabet, then every permutation $\sigma\in\SG(\alphA)$ gives rise
to an automorphism
on the free algebra $\K\<\alphA\>$. Two words $u=u_1u_2\cdots u_n$ and $v=v_1v_2\cdots v_n$
are in the same orbit whenever $u_i=u_j \Leftrightarrow v_i=v_j$, i.e.,
$\ker u=\ker v$.
Thus, the
orbits are in one-to-one correspondence with set partitions into at most
$|\alphA|$ blocks.

\subsubsection{Monomial basis}
Assuming as above
that $\alphA$ is infinite, we obtain an algebra based on all set partitions,
defining  the monomial basis by  symmetrization
\marginpar{$\m_\pi$}
\begin{equation}
  \label{eq:def:WSym:mpi}
  \m_\pi(\alphA) 
  = \sum_{\ker w=\pi} w
  .
\end{equation}
For example,  the word $a_1a_2a_1a_2^2$ has kernel $\pi=13/245=
\begin{picture}(26,6.5)(1,0)
\put(2,0){\line(0,1){3.5}}
\put(8,0){\line(0,1){6.5}}
\put(14,0){\line(0,1){3.5}}
\put(20,0){\line(0,1){6.5}}
\put(26,0){\line(0,1){6.5}}
\put(2,3.5){\line(1,0){12}}
\put(8,6.5){\line(1,0){18}}
\end{picture}
$ and its
symmetrization yields the element
$$
\m_{
  \scriptpic{\begin{picture}(26,6.5)(1,0)
\put(2,0){\line(0,1){3.5}}
\put(8,0){\line(0,1){6.5}}
\put(14,0){\line(0,1){3.5}}
\put(20,0){\line(0,1){6.5}}
\put(26,0){\line(0,1){6.5}}
\put(2,3.5){\line(1,0){12}}
\put(8,6.5){\line(1,0){18}}
\end{picture}}
} = \sum_{i\ne j} a_ia_ja_ia_j^2
.
$$
The abelianization $\theta(\m_\pi)$ yields the element $m_\pi$ from
Doubilet's formulas   \eqref{eq:doubilet},
which is a multiple of the commutative
monomial symmetric function $m_{\lambda(\pi)}$.

\subsubsection{Coproduct}
One can introduce a bialgebra structure by means of the coproduct
\begin{equation}
\Delta f(\alphA) = f(\alphA'+\alphA'')
\end{equation}
where $\alphA'+\alphA''$ denotes the disjoint union of two mutually commuting copies of
$\alphA$, i.e., $a'=a\otimes 1$ and $a''=1\otimes a$.
This coproduct is obviously cocommutative and using the standardization of
partitions introduced in §\ref{ssec:setpartitions} the coproduct of a
monomial function can be written as
\begin{equation}
  \Delta \m_\pi
  = \sum_{\pi'\dot{\cup}\pi''=\pi}\m_{\std(\pi')}\otimes \m_{\std(\pi'')}
  .
\end{equation}
For example,
\begin{multline*}
\m_{\scriptpic{\begin{picture}(38,6.5)(1,0)
\put(2,0){\line(0,1){3.5}}
\put(8,0){\line(0,1){6.5}}
\put(14,0){\line(0,1){3.5}}
\put(20,0){\line(0,1){6.5}}
\put(26,0){\line(0,1){3.5}}
\put(32,0){\line(0,1){6.5}}
\put(38,0){\line(0,1){3.5}}
\put(2,3.5){\line(1,0){12}}
\put(8,6.5){\line(1,0){24}}
\put(26,3.5){\line(1,0){12}}
\end{picture}}}
=
\m_{\scriptpic{
    \begin{picture}(38,6.5)(1,0)
\put(2,0){\line(0,1){3.5}}
\put(8,0){\line(0,1){6.5}}
\put(14,0){\line(0,1){3.5}}
\put(20,0){\line(0,1){6.5}}
\put(26,0){\line(0,1){3.5}}
\put(32,0){\line(0,1){6.5}}
\put(38,0){\line(0,1){3.5}}
\put(2,3.5){\line(1,0){12}}
\put(8,6.5){\line(1,0){24}}
\put(26,3.5){\line(1,0){12}}
\end{picture}}}
\otimes
1
+
\m_{\scriptpic{\begin{picture}(26,6.5)(1,0)
\put(2,0){\line(0,1){3.5}}
\put(8,0){\line(0,1){3.5}}
\put(14,0){\line(0,1){6.5}}
\put(20,0){\line(0,1){3.5}}
\put(26,0){\line(0,1){6.5}}
\put(2,3.5){\line(1,0){18}}
\put(14,6.5){\line(1,0){12}}
\end{picture}}}
\otimes
\m_{\scriptpic{
    \begin{picture}(8,3.5)(1,0)
\put(2,0){\line(0,1){3.5}}
\put(8,0){\line(0,1){3.5}}
\put(2,3.5){\line(1,0){6}}
\end{picture}}}
+
\m_{\scriptpic{
    \begin{picture}(20,3.5)(1,0)
\put(2,0){\line(0,1){3.5}}
\put(8,0){\line(0,1){3.5}}
\put(14,0){\line(0,1){3.5}}
\put(20,0){\line(0,1){3.5}}
\put(2,3.5){\line(1,0){6}}
\put(14,3.5){\line(1,0){6}}
\end{picture}}}
\otimes
\m_{\scriptpic{
    \begin{picture}(14,3.5)(1,0)
\put(2,0){\line(0,1){3.5}}
\put(8,0){\line(0,1){3.5}}
\put(14,0){\line(0,1){3.5}}
\put(2,3.5){\line(1,0){12}}
\end{picture}}}
+
\m_{\scriptpic{\begin{picture}(26,6.5)(1,0)
\put(2,0){\line(0,1){3.5}}
\put(8,0){\line(0,1){6.5}}
\put(14,0){\line(0,1){3.5}}
\put(20,0){\line(0,1){6.5}}
\put(26,0){\line(0,1){6.5}}
\put(2,3.5){\line(1,0){12}}
\put(8,6.5){\line(1,0){18}}
\end{picture}}}
\otimes
\m_{\scriptpic{\begin{picture}(8,3.5)(1,0)
\put(2,0){\line(0,1){3.5}}
\put(8,0){\line(0,1){3.5}}
\put(2,3.5){\line(1,0){6}}
\end{picture}}}
\\
+
\m_{\scriptpic{
    \begin{picture}(8,3.5)(1,0)
\put(2,0){\line(0,1){3.5}}
\put(8,0){\line(0,1){3.5}}
\put(2,3.5){\line(1,0){6}}
\end{picture}}}
\otimes
\m_{\scriptpic{\begin{picture}(26,6.5)(1,0)
\put(2,0){\line(0,1){3.5}}
\put(8,0){\line(0,1){6.5}}
\put(14,0){\line(0,1){3.5}}
\put(20,0){\line(0,1){6.5}}
\put(26,0){\line(0,1){6.5}}
\put(2,3.5){\line(1,0){12}}
\put(8,6.5){\line(1,0){18}}
\end{picture}}}
+
\m_{\scriptpic{\begin{picture}(14,3.5)(1,0)
\put(2,0){\line(0,1){3.5}}
\put(8,0){\line(0,1){3.5}}
\put(14,0){\line(0,1){3.5}}
\put(2,3.5){\line(1,0){12}}
\end{picture}}}
\otimes
\m_{\scriptpic{\begin{picture}(20,3.5)(1,0)
\put(2,0){\line(0,1){3.5}}
\put(8,0){\line(0,1){3.5}}
\put(14,0){\line(0,1){3.5}}
\put(20,0){\line(0,1){3.5}}
\put(2,3.5){\line(1,0){6}}
\put(14,3.5){\line(1,0){6}}
\end{picture}}}
+
\m_{\scriptpic{\begin{picture}(8,3.5)(1,0)
\put(2,0){\line(0,1){3.5}}
\put(8,0){\line(0,1){3.5}}
\put(2,3.5){\line(1,0){6}}
\end{picture}}}
\otimes
\m_{\scriptpic{\begin{picture}(26,6.5)(1,0)
\put(2,0){\line(0,1){3.5}}
\put(8,0){\line(0,1){3.5}}
\put(14,0){\line(0,1){6.5}}
\put(20,0){\line(0,1){3.5}}
\put(26,0){\line(0,1){6.5}}
\put(2,3.5){\line(1,0){18}}
\put(14,6.5){\line(1,0){12}}
\end{picture}}}
+
1\otimes
\m_{\scriptpic{\begin{picture}(38,6.5)(1,0)
\put(2,0){\line(0,1){3.5}}
\put(8,0){\line(0,1){6.5}}
\put(14,0){\line(0,1){3.5}}
\put(20,0){\line(0,1){6.5}}
\put(26,0){\line(0,1){3.5}}
\put(32,0){\line(0,1){6.5}}
\put(38,0){\line(0,1){3.5}}
\put(2,3.5){\line(1,0){12}}
\put(8,6.5){\line(1,0){24}}
\put(26,3.5){\line(1,0){12}}
\end{picture}}}
\end{multline*}

\subsubsection{Product}
On the other hand the product is more complicated in terms of set partitions
$$
\m_\pi\m_\rho = \sum_{\sigma\land(\hat{1}_m|\hat{1}_n) = \pi|\rho} \m_\sigma
.
$$
For example,
$$
\m_{\scriptpic{
    \begin{picture}(14,3.5)(1,0)
\put(2,0){\line(0,1){3.5}}
\put(8,0){\line(0,1){3.5}}
\put(14,0){\line(0,1){3.5}}
\put(2,3.5){\line(1,0){0}}
\put(8,3.5){\line(1,0){6}}
\end{picture}}}
\m_{\scriptpic{
    \begin{picture}(14,6.5)(1,0)
\put(2,0){\line(0,1){6.5}}
\put(8,0){\line(0,1){3.5}}
\put(14,0){\line(0,1){6.5}}
\put(8,3.5){\line(1,0){0}}
\put(2,6.5){\line(1,0){12}}
\end{picture}}}
=
  \m_{\scriptpic{\begin{picture}(32,6.5)(1,0)
\put(2,0){\line(0,1){6.5}}
\put(8,0){\line(0,1){3.5}}
\put(14,0){\line(0,1){3.5}}
\put(20,0){\line(0,1){6.5}}
\put(26,0){\line(0,1){3.5}}
\put(32,0){\line(0,1){6.5}}
\put(8,3.5){\line(1,0){18}}
\put(2,6.5){\line(1,0){30}}
\end{picture}}} 
+
\m_{\scriptpic{\begin{picture}(32,6.5)(1,0)
\put(2,0){\line(0,1){6.5}}
\put(8,0){\line(0,1){3.5}}
\put(14,0){\line(0,1){3.5}}
\put(20,0){\line(0,1){6.5}}
\put(26,0){\line(0,1){3.5}}
\put(32,0){\line(0,1){6.5}}
\put(8,3.5){\line(1,0){6}}
\put(26,3.5){\line(1,0){0}}
\put(2,6.5){\line(1,0){30}}
\end{picture}}} 
+
\m_{\scriptpic{\begin{picture}(32,6.5)(1,0)
\put(2,0){\line(0,1){3.5}}
\put(8,0){\line(0,1){6.5}}
\put(14,0){\line(0,1){6.5}}
\put(20,0){\line(0,1){6.5}}
\put(26,0){\line(0,1){3.5}}
\put(32,0){\line(0,1){6.5}}
\put(2,3.5){\line(1,0){24}}
\put(8,6.5){\line(1,0){24}}
\end{picture}}} 
+
\m_{\scriptpic{\begin{picture}(32,6.5)(1,0)
\put(2,0){\line(0,1){3.5}}
\put(8,0){\line(0,1){6.5}}
\put(14,0){\line(0,1){6.5}}
\put(20,0){\line(0,1){6.5}}
\put(26,0){\line(0,1){3.5}}
\put(32,0){\line(0,1){6.5}}
\put(2,3.5){\line(1,0){0}}
\put(26,3.5){\line(1,0){0}}
\put(8,6.5){\line(1,0){24}}
\end{picture}}} 
+
\m_{\scriptpic{\begin{picture}(32,6.5)(1,0)
\put(2,0){\line(0,1){6.5}}
\put(8,0){\line(0,1){3.5}}
\put(14,0){\line(0,1){3.5}}
\put(20,0){\line(0,1){3.5}}
\put(26,0){\line(0,1){6.5}}
\put(32,0){\line(0,1){3.5}}
\put(8,3.5){\line(1,0){6}}
\put(2,6.5){\line(1,0){24}}
\put(20,3.5){\line(1,0){12}}
\end{picture}}} 
+
\m_{\scriptpic{\begin{picture}(32,6.5)(1,0)
\put(2,0){\line(0,1){3.5}}
\put(8,0){\line(0,1){3.5}}
\put(14,0){\line(0,1){3.5}}
\put(20,0){\line(0,1){6.5}}
\put(26,0){\line(0,1){3.5}}
\put(32,0){\line(0,1){6.5}}
\put(2,3.5){\line(1,0){0}}
\put(8,3.5){\line(1,0){18}}
\put(20,6.5){\line(1,0){12}}
\end{picture}}} 
+
\m_{\scriptpic{\begin{picture}(32,6.5)(1,0)
\put(2,0){\line(0,1){3.5}}
\put(8,0){\line(0,1){3.5}}
\put(14,0){\line(0,1){3.5}}
\put(20,0){\line(0,1){6.5}}
\put(26,0){\line(0,1){3.5}}
\put(32,0){\line(0,1){6.5}}
\put(2,3.5){\line(1,0){0}}
\put(8,3.5){\line(1,0){6}}
\put(26,3.5){\line(1,0){0}}
\put(20,6.5){\line(1,0){12}}
\end{picture}}}
.
$$
It becomes simpler when we switch to the noncommutative power sums.
Indeed the basis of noncommutative power sums $\Phi^\pi$, defined by sums over
initial intervals
\marginpar{$\Phi^\pi$}
\begin{equation}
  \label{eq:WSym:Phi=summ}
  \Phi^\pi =\sum_{\pi'\SFleq\pi}\m_{\pi'}
  ,
\end{equation}
where $\SFleq$ is the  refinement order   \eqref{eq:SFleq},
is multiplicative with respect to concatenation of set partitions:
$$
\Phi^\pi\Phi^\rho = \Phi^{\pi|\rho}
$$
and hence $\WSym$ is freely generated by the elements $\Phi^\pi$ such that
$\pi$ is irreducible, i.e.,
the finest interval partition dominating $\pi$ has only one block.

\subsubsection{Embedding into $\WQSym$}
The algebra $\WSym$ is a subalgebra of $\WQSym$ via the embedding
\begin{equation}
  \label{eq:WSyminWQSym}
\mathbf{m}_\pi \mapsto \sum_{\substack{u\\ \ker u=\pi }} \mathbf{M}_u  = \sum_\sigma \mathbf{M}_{\sigma(w)}
\end{equation}
where $u$ runs over all packed words,
and $\sigma$ runs over all permutations of the letters of the
restricted growth function $w=\rgf(\pi)$ which encodes the partition $\pi$.
For example, the monomial corresponding to the partition $\pi=
\begin{picture}(26,6.5)(1,0)
\put(2,0){\line(0,1){3.5}}
\put(8,0){\line(0,1){6.5}}
\put(14,0){\line(0,1){3.5}}
\put(20,0){\line(0,1){3.5}}
\put(26,0){\line(0,1){6.5}}
\put(2,3.5){\line(1,0){12}}
\put(20,3.5){\line(1,0){0}}
\put(8,6.5){\line(1,0){18}}
\end{picture}$ is mapped to
$$
\mathbf{m}_{\scriptpic{\begin{picture}(26,6.5)(1,0)
\put(2,0){\line(0,1){3.5}}
\put(8,0){\line(0,1){6.5}}
\put(14,0){\line(0,1){3.5}}
\put(20,0){\line(0,1){3.5}}
\put(26,0){\line(0,1){6.5}}
\put(2,3.5){\line(1,0){12}}
\put(20,3.5){\line(1,0){0}}
\put(8,6.5){\line(1,0){18}}
\end{picture}}}
\mapsto
\mathbf{M}_{12132}
+
\mathbf{M}_{13123}
+
\mathbf{M}_{21231}
+
\mathbf{M}_{23213}
+
\mathbf{M}_{31321}
+
\mathbf{M}_{32312}
.
$$

\subsubsection{Dual}
\label{sssec:WSym*}
The graded dual of $\WSym$ is a commutative
algebra and thus $\WSym^*\not\simeq\WSym$.
It is isomorphic to the algebra  $\PiQSym$ defined in 
\cite[§3.5.1]{HNT08}.
\marginpar{$N_\pi$}
Let $N_\pi$
be the dual basis of $\m_\pi$ in the (commutative) graded  dual
$\WSym^*\simeq \PiQSym$
and let $\phi_\pi$
\marginpar{$\phi_\pi$}
be the dual basis of $\Phi^\pi$.
Then,
\begin{equation}
  \label{eq:Nphi}
N_{\pi'} =\sum_{\pi\SFgeq\pi'}\phi_{\pi}.
\end{equation}

Dually to the embedding  \eqref{eq:WSyminWQSym},
$\WSym^*$ is a quotient of $\WQSym^*$,
with respect to the relation on packed words defined by
$\N_u\equiv \N_v$ $\iff$ $\ker u=\ker v$, i.e.,
$u=\sigma(v)$ for some permutation
$\sigma$ of the alphabet of $u$ (e.g.,
$
\N_{12132}
\equiv
\N_{13123}
\equiv
\N_{21231}
\equiv
\N_{23213}
\equiv
\N_{31321}
\equiv
\N_{32312}
$)
and the equivalence classes are parametrized by
the set partitions corresponding to any of the set compositions encoded by equivalent packed words.
\subsubsection{Internal product}
The internal product of $\WQSym^*$ descends to an internal product on $\WSym^*$, which is given by the meet
operation in the lattice of set partitions:
\begin{equation}
  \label{eq:Npi*Ntau=Npiwedgetau}
N_\pi*N_\tau = N_{\pi\wedge\tau}
\end{equation} 
where $\pi\wedge\tau$ is the coarsest partition which is finer than $\pi$ and $\tau$.
Endowed with this product, the homogenous component $\WSym_n^*$ is isomorphic to the Möbius algebra
of the partition lattice $\Pi_n$,
see §\ref{sssec:moebiusalg}.
It has been shown by Solomon
\cite{Sol76}
that  $\{\phi_\pi|\pi\in\Pi_n\}$
is a complete set of orthogonal idempotents of $\WSym^*_n$.

%

As a consequence, if $\pi$ is not the trivial partition $\{12\ldots n\}$,
\begin{equation}
  \label{eq:prodNphi}
N_\pi*\phi_{\{12\cdots n\}}=0.
\end{equation}

\subsection{Noncommutative version of Good's formula}

Although the identity
\eqref{eq:omegacom}
(resp.~\eqref{eq:exch:Kn=sumphimu})
is essentially trivial, it has a not-so-trivial noncommutative
analogue for noncommutative symmetric functions. It is proved in \cite[Prop. 8.6]{NCSF2}
that
\begin{equation}
S_n(\Omega_n \alphA) = K_n(\omega)
\end{equation}
where $K_n(\omega)$ is Klyachko's element, a famous Lie (quasi-)idempotent.
This fact leads to an interesting interpretation of
the formula of Good and Cartier
    \eqref{eq:omegacom}.
To understand it, we have to follow the chain of morphisms
\begin{equation}
\Sym\hookrightarrow \WQSym^* \twoheadrightarrow \WSym^* \simeq \PiQSym.
\end{equation}

The first embedding $i$ is given by
  \eqref{eq:embed:Symn->WQSym*}.
The projection onto $\WSym^*$ is dual to the
inclusion
  \eqref{eq:WSyminWQSym}
  of $\WSym$ into $\WQSym$
  and therefore given by
  \begin{equation}
    \label{eq:p:WQsumtoWSym}
  \begin{aligned}
    p:\WQSym^*&\to\WSym\\
      \N_u\mapsto N_\pi
  \end{aligned}    
  \end{equation}
 where $N_\pi$ is dual to the
 monomial basis of $\WSym$, and $\pi=\ker u$ is the support
 of 
 the set composition encoded by the packed word $u$.

In \cite{Leh04},  cumulants for an exchangeability system are defined by a
noncommutative version of the formula of Good and Cartier \eqref{eq:omegacom}.
We shall see that they can
be interpreted as the image of $(n-1)!K_n(\omega)$ under the composition
$\xi=p\circ i$ of
the maps     \eqref{eq:p:WQsumtoWSym} and    \eqref{eq:embed:Symn->WQSym*}.

Since $\xi$ is valued in a commutative algebra, it factors through $\sym(\alphX)$, and all noncommutative
power sums have the same image. The choice of Klyachko's element is therefore arbitrary. 

It remains to determine its image. Let $\Phi^\pi = \sum_{\pi\SFleq\tau}\M_\tau$ be the power-sum basis
of $\WSym$, and let $\phi_\pi$ be its dual basis. Under the abelianization map 
\begin{equation}
  \label{eq:WSym:abelianization}
\theta:\ \WSym(\alphA) \twoheadrightarrow \sym(\alphX)
\end{equation}
mapping $a_i\mapsto X_i$ and $f(\alphA)\mapsto f(\alphX)$, we have
\begin{equation}
\M_\pi(\alphA) \mapsto \prod_i r_i(\lambda)! \cdot m_\lambda(\alphX) 
\end{equation}  
where $\lambda$ is the integer partition associated with $\pi$,
$m_\lambda(\alphX)$ is the monomial symmetric function,  $r_i( \lambda)$ is the
multiplicity of the part $i$ in $\lambda$, and
\begin{equation}
\Phi^\pi(\alphA) \mapsto  p_\lambda(\alphX). 
\end{equation}
Dually,
\begin{equation}
  \label{eq:theta*hlambda}
\theta^*(h_\lambda) = \prod_i r_i(\lambda)! \sum_{\Lambda(\pi)=\lambda}N_\pi,
\end{equation}
and
\begin{equation}
\theta^*(p_\lambda) = z_\lambda \sum_{\Lambda(\pi)=\lambda}\phi_\pi.
\end{equation}
The image of  any Lie idempotent of $\Sym$ in $\WSym^*$ is therefore the same as that
of $\varphi_n = \frac1n\Phi_n$, which is
\begin{equation}
\varphi_n = 
\sum_{I\vDash n} \frac{(-1)^{\ell(I)-1}}{\ell(I)}S^I
\mapsto    \sum_{|u|=n}\frac{(-1)^{\max(u)-1}}{\max(u)}\N_u
\mapsto 
\sum_{\pi\vdash [n]}\frac{(-1)^{\ell(\pi)-1}}{\ell(\pi)}\ell(\pi)!N_\pi.
\end{equation}
Since this must also be equal to
\begin{equation} 
\theta^*(\frac1n\, p_n) =\phi_{\{1,2,\ldots, n \} }, 
\end{equation}
and since, by definition of the Möbius function of the lattice of partitions
\begin{equation}
  \label{eq:WSym*:phi=sumNmu}
\phi_{\{1,2,\ldots, n \} }=\sum_{\pi\vdash [n]}N_\pi\mu(\pi,\hat 1_n)
\end{equation}
we recover the classical fact that
\begin{equation}
  \label{eq:moebiusmu}
\mu(\pi,\hat 1_n)=(-1)^{\ell(\pi)-1}(\ell(\pi)-1)!
\end{equation}
by merely contemplating a chain of morphisms.

\subsection{Cumulants for exchangeability systems}

An exchangeability system together with a sequence $\mathcal{X}=(X_i)_{i\ge 1}$ determines
a linear form $\hat\varphi_{\mathcal{X}}$ on $\WSym^*$ 
\begin{equation}
\hat\varphi_{\mathcal{X}}(N_\pi) := \varphi(X_1^{(u_1)}X_2^{(u_2)}\dotsm X_n^{(u_n)}).
\end{equation}
where $u$ is any packed word such that $\ker u=\pi$,
see §\ref{sssec:exchangeability}.

The partitioned moments are
\begin{equation}
\varphi_\pi(X_1X_2\dotsm X_n) = \hat\varphi_{\mathcal{X}}(N_\pi),
\end{equation}
and the cumulants defined in  
  \eqref{eq:phiNX=NKnX}
  are
\begin{equation}
K_\pi(X_1,X_2,\ldots,X_n) =\hat\varphi_{\mathcal{X}}(\phi_\pi).
\end{equation}
By \eqref{eq:Nphi}, we have \cite[Prop. 2.7]{Leh04}
\begin{equation}
  \label{eq:phipi=sumksig}
\varphi_\pi(X_1,X_2,\ldots,X_n)=\sum_{\pi \SPleq\sigma}K_\sigma(X_1,X_2,\ldots,X_n).
\end{equation}
In the case of exchangeability systems the  independence condition  defined in \eqref{eq:Sindep} becomes
\begin{equation} %
  \label{eq:Exch:phiNpi=ph}
\hat\varphi_{\mathcal{X}}(N_\pi)=\hat\varphi_{\mathcal{X}}(N_\pi*N_\sigma)=\hat\varphi_{\mathcal{X}}(N_{\pi\wedge\sigma}).
\end{equation}
This implies the vanishing of mixed cumulants:
if $\pi$ is a partition if $[n]$
into two nontrivial blocks $b_1,b_2$ such that $\{X_i|i\in b_1\}$ and $\{X_i|i\in b_2\}$
are independent, then the cumulant 
\begin{equation}
K_n(X_1,X_2,\ldots,X_n)=\hat\varphi_{\mathcal{X}}(\phi_{1^n})=\hat\varphi_{\mathcal{X}}(N_\pi*\phi_{1^n}) = 0
\end{equation}
by \eqref{eq:prodNphi}.

\subsection{Cumulants for spreadability systems}
Monotone independence is asymmetric and does not fit into the frame of
${\mathcal E}$-independence, which was one of the motivations to introduce the
notion of spreadability and ${\mathcal S}$-independence.
In this case, the requirement that mixed cumulants vanish is too
strong, and is replaced by the axiom of extensivity (K1$^\prime$)
introduced in §\ref{sec:cumulants}, i.e., 
the condition that they should be eigenfunctions of Rota's dot
multiplication by an integer:
\begin{equation}
K_u(N.X_1,N.X_2,\ldots,N.X_n)=N^{\max(u)}K_u(X_1,X_2,\ldots,X_n)
\end{equation}
where $u$ is a packed word, and $N.X=X^{(1)}+X^{(2)}+\cdots+X^{(N)}$ is the sum of
$N$ independent copies of $X$. 
We shall see that the dot operation translates as multiplication of the alphabet by $N$
in the relevant Hopf algebras. The operator $f(\alphA)\mapsto f(N\alphA)$
from   \eqref{eq:def:f(kA)}
is semisimple, and its spectral
projectors are know as the \emph{Eulerian idempotents}.

\section{Review of the Eulerian algebra}

\subsection{Basics}
The Eulerian algebra is a commutative subalgebra of dimension $n$ of the group algebra of
the symmetric group $\SG_n$, and in fact of its descent algebra $\Sigma_n$. It was apparently
first  introduced in \cite{BMP} under the name {\em algebra of permutors}\footnote{
A self-contained and elementary presentation of the main results of \cite{BMP}
can be found in \cite{NCSF1}.}.
It is spanned by the Eulerian idempotents, or, equivalently, by the sums of
permutations having the same number of descents \cite{Reu86}. The name ``Eulerian idempotents''
is apparently due to Loday \cite{Lod}. The representations associated with these idempotents have been first investigated
by Reutenauer \cite{Reu86}. 

It is easier to work with all symmetric groups at the same time, with the
help of generating functions. Recall from
§\ref{sssec:SYM:internalproduct}
that the algebra of noncommutative
symmetric functions $\Sym$ is endowed with an internal product $*$,
for which each homogeneous component $\Sym_n$ is anti-isomorphic to $\Sigma_n$
\cite[§5.1]{NCSF1} and has unit element $S_n$.

The {\em Eulerian idempotents} $E_n^{[k]}$ are the homogenous components
of degree $n$ in the series $E^{[k]}$ defined by
\begin{equation}
  \label{eq:def:E[k]}
\sigma_t(\alphA)^x=\sum_{k\ge 0}x^k E^{[k]}(\alphA),
\end{equation}
(see \cite[§5.3]{NCSF1});
that is, $E_n^{[k]}$
\marginpar{$E_n^{[k]}$}
is the coefficient of $t^n$ in the formal expansion
of
$$
\frac{1}{k!}\log (1+\sum_{m=1}^\infty S_mt^m)
.
$$
Then
\begin{equation}
E_n^{[k]}*E_n^{[l]} = \delta_{kl}E_n^{[k]}\,,\quad\text{and}\quad\sum_{k=1}^nE_n^{[k]}=S_n,
\end{equation}
so that the $E_n^{[k]}$ form an orthogonal partition of unity and span a commutative $n$-dimensional $*$-subalgebra of
$\Sym_n$, denoted by $\Eul_n$ and called the \emph{Eulerian subalgebra}.

\subsection{Noncommutative Eulerian polynomials}

The {\em noncommutative Eulerian polynomials} are defined by \cite[§5.4]{NCSF1}
\begin{equation}
{\mathcal A}_n(t) =
\ \sum_{k=1}^n \ t^k\, \Big(
\sum_{{\scriptstyle |I|=n}\atop{\scriptstyle \ell(I)=k}} R_I \,
\Big)
=
\ \sum_{k=1}^n \ {\bf A}(n,k)\, t^k \, ,
\end{equation}
where $R_I$
\marginpar{$R_I$}
is the ribbon basis \cite[§3.2]{NCSF1}
$$
R_I = \sum_{\Desc(J)\subseteq \Desc(I)} (-1)^{\ell(I)-\ell(J)} S^J
.
$$

The following facts can be found (up to a few misprints\footnote{
Eqs. (93) and (97) of \cite{NCSF1} should be read as (\ref{A2S}) and (\ref{Wor}) of the present paper.
}) in \cite{NCSF1}.
The generating series of the ${\mathcal A}_n(t)$ is
\begin{equation}
{\mathcal A}(t) := \ \sum_{n\ge 0} \, {\mathcal A}_n(t)
=
(1-t) \, \left( 1 - t\, \sigma_{1-t} \right)^{-1} \,,
\end{equation}
where $\sigma_{1-t}=\sum_{k=0}^\infty (1-t)^kS_k$ is the generating function   \eqref{eq:def:sigmat}.

\medskip
Let ${\mathcal A}_n^*(t) = (1-t)^{-n}\, {\mathcal A}_n(t)$.
Then,
\begin{equation}
{\mathcal A}^*(t)
:=
\ \sum_{n\ge 0} \, {\mathcal A}_n^*(t)
=
\sum_{I} \
\left( \displaystyle {t \over 1-t} \right)^{\ell(I)} \, S^I \ .
\end{equation}
This last formula can also be written in the form
\begin{equation} \label{GEN*}
{\mathcal A}^*(t)
=
\ \sum_{k\ge 0} \ \left(
{t\over 1-t}\right)^k \left( S_1+S_2+S_3+\cdots\, \right)^k
\end{equation}
or
\begin{equation}\label{GEN_A}
{1\over 1-t\, \sigma_1(\alphA)}
=
\ \sum_{n\ge 0}\ {{\mathcal A}_n(t)\over (1-t)^{n+1}} \ .
\end{equation}
Let $S^{[k]}=\sigma_1(\alphA)^k$ be the coefficient of $t^k$ in this series. In
degree $n$,
\marginpar{$S_n^{[k]}$}
\begin{equation}\label{S2E}
S_n^{[k]}=\sum_{I\vDash n, \ell(I)\le k}{k\choose \ell(I)}S^I=\sum_{i=1}^nk^iE_n^{[i]}\,.
\end{equation}
This is another basis of $\Eul_n$. 
Expanding the factors $(1-t)^{-(n+1)}$ on the right-hand side
of (\ref{GEN_A}) by the binomial
theorem, and taking the coefficient of $t^k$ in the term of
weight $n$ on both sides, we get
\begin{equation}
S_n^{[k]}
=
\ \sum_{i=0}^k \ {n+i \choose i} \, {\bf A}(n,k-i) \, .
\end{equation}
Conversely,
\begin{equation}
{ {\mathcal A}_n(t) \over (1-t)^{n+1} }
=
\ \sum_{k\ge 0}\ t^k\, S_n^{[k]} \ ,
\end{equation}
so that
\begin{equation}\label{A2S}
{\bf A}(n,p)
=
\ \sum_{i=0}^p\ (-1)^i\, {n+1\choose i}\, S_n^{[p-i]} \ .
\end{equation}
The expansion of the $E_n^{[k]}$
in the basis ${\bf A}(n,i)$, which is a noncommutative analog
of Worpitzky's identity (see
\cite[§2.6]{FoataSchuetzenberger:theorie:1970}
or
\cite[(6.37)]{GrahamKnuthPatashnik:2:1994}) is
\begin{equation}\label{Wor}
\sum_{k=1}^n \ x^k\, E_n^{[k]}
=
\ \sum_{i=1}^n \ {x+n-i\choose n}\, {\bf A}(n,i) \ .
\end{equation}
Indeed, when $x$ is a positive integer $N$,
\begin{equation}
\sum_{k=1}^n \ N^k\, E_n^{[k]}
= S_n(N\alphA) = \sum_{I\vDash n}F_I(N)R_I(\alphA)
\end{equation} 
where $F_I$
\marginpar{$F_I$}
are the \emph{fundamental quasi-symmetric functions},
and for a composition $I=(i_1,\ldots,i_r)$ of $n$,
\begin{equation}
F_I(N)={N+n-r\choose n}\,.
\end{equation}


\section{Adams operations of $\Sym$ and their substitutes}

On any bialgebra $\H$ with multiplication $\mu$ and comultiplication $\Delta$, one
can define the {\it Adams operations} \cite[§4.1]{CartierPatras:2021}
\marginpar{$\Psi^k$}
\begin{equation}
\Psi^k(x) = \mu_k\circ \Delta^k(x)
\end{equation}
where $\Delta^k$ is the iterated coproduct with values in $\H^{\otimes k}$ and $\mu_k$
the multiplication map $\H^{\otimes k}\rightarrow \H$. In other terms, $\Psi^k$ is the
$k$th convolution power of the identity. 

On ordinary symmetric functions, these
operations act by multiplication by $k$ of the alphabet, and are therefore algebra
morphisms. However, on noncommutative symmetric functions, the $\Psi^k$ are not anymore
multiplicative (e.g., $\Psi^2(S^{1,2})\ne\Psi^2(S^{1})\Psi^2(S^{2}) $), and therefore of minor interest.
One can replace them by the algebra morphisms $\psi_k:\ f(\alphA)\mapsto
f(k\alphA)$ with
\marginpar{$\psi_k$, $f(k\alphA)$}
\begin{equation}
  \label{eq:def:f(kA)}
  f(k\alf{A}) = f(\alf{A}'+\alf{A}''+\dots+\alf{A}^{(k)})
\end{equation}
where $\alf{A}'$, $\alf{A}''$,\ldots, $\alf{A}^{(k)}$ are mutually commuting copies
of the alphabet $\alf{A}$.
For example, from   \eqref{eq:Sym:sigma(A+B)} we infer that
%
\begin{equation} 
\sigma_1(k\alphA) = \sigma_1(\alphA)^k.
\end{equation}
We will see that the operators $\psi_k$  are diagonalized by the Eulerian idempotents. 
Indeed
 for any $f\in\Sym$ the identity
 \begin{equation}
   \label{eq:f(kA)=f(A)*sigma1}
f(k\alphA) = f(\alphA) * \sigma_1(k\alphA)
\end{equation}
holds.
To see this, it suffices to show two things:
that the right hand side of    \eqref{eq:f(kA)=f(A)*sigma1} defines
an algebra morphism and that the identity holds whenever $f$ is substituted
with a  generator $S_n$, or equivalently, with their direct sum $\sigma_1(\alphA)$.

For the first requirement, observe that  $\sigma_1(\alphA)$, being the
direct sum of the units for the internal product on $\Sym_n$, satisfies
$f*\sigma_1(\alphA)=f$. 
On the other hand, $\sigma_1(k\alphA)$ is group-like
and therefore from the splitting formula \eqref{eq:Sym:internal-coprod}
we infer that
\begin{align*}
  (f\cdot g)*\sigma_1(k\alphA)
  &= \mu_2((f\otimes g)*_2(\sigma_1(k\alphA)\otimes\sigma_1(k\alphA)))
    \\
  &=   (f*\sigma_1(k\alphA))\cdot  (g*\sigma_1(k\alphA))
    .
\end{align*}

The second requirement follows trivially from  \eqref{eq:sigma1*f=f}
with $f=\sigma_1(k\alphA)$.

  Finally with   \eqref{eq:def:E[k]} we conclude that
  $$
f(k\alphA)
 = f(\alphA)*\sum_{r\ge 1} k^r E^{[r]}.
 $$


Thus, a simultaneous eigenbasis of the $\psi_k$ is for example
\begin{equation}
K_I := S^I * E^{[\ell(I)]}
\end{equation} 
which satisfy
\begin{equation}
K_I(k\alphA) = k^{\ell(I)}K_I(\alphA). 
\end{equation}
The basis $K_I$ is actually multiplicative: if $I=(i_1,\ldots,i_r)$, then by
the splitting formula \eqref{eq:Sym:internal-coprod}
\begin{align}\label{eq:mult}
K_I& = [k^r] S^I * \sigma_1(\alphA)^k = [k^r]\mu_r[(S_{i_1}\otimes\cdots\otimes S_{i_r})*_r (\sigma_1(k\alphA)\otimes\cdots\otimes \sigma_1(k\alphA)]\\
&= (S_{i_1}*E^{[1]})\cdots  (S_{i_r}*E^{[1]}) = K_{i_1}\cdots K_{i_r}.
\end{align} 
Of course, $S_i*E_i^{[1]}=E_i^{[1]}=\frac1i \Phi_i$, so that $K_I$ is just a scaled version of
the classical basis $\Phi^I$.

Since
\begin{equation} 
\sigma_1(k\alphA) = [1+(S_1+S_2+S_3+\cdots)]^k = \sum_{I}{k\choose \ell(I)}S^I,
\end{equation}
we have the simple closed form
\begin{equation}
S^I(k\alphA) = \sum_{J\ge I}\beta_k(J,I) S^J,
\end{equation} 
with
\begin{equation}
\beta_k(J,I) := \prod_p{k\choose \ell(J_p)},
\end{equation}
where $J=(J_1J_2\cdots J_r)$ is a concatenation of compositions $J_p$ of weight $i_p$ (by definition
of the refinement order).

Applying \cite[Prop. 4.9]{NCSF1}, we get the expression
\begin{equation}
K_I=\sum_{J\ge I}\frac{(-1)^{\ell(J)-\ell(I)}}{\ell(J,I)}S^J,
\end{equation}
where $\ell(J,I)=\prod_{p=1}^r\ell(J_p)$.

\section{Extension to $\WQSym^*$}

The identification of $\Sym_n$ with the (opposite) descent algebra of $\SG_n$ can be refined
as follows. 
We have seen that $\WQSym^*$ can be identified with the (opposite) Solomon-Tits algebra. The dual basis
$\N_u=\M_u^*$ of the monomial basis of $\WQSym$ has the internal product rule
  \eqref{eq:WQSym*:internalproduct}
and the embedding of the Solomon algebra into the Solomon-Tits algebra is given by
the Hopf embedding
  \eqref{eq:embed:Symn->WQSym*}
of $\Sym$ into $\WQSym^*$.
In particular, $S_n(k\alphA)$ and the Eulerian idempotents can be interpreted as elements of $\WQSym_n^*$, and
one can define a new basis
\marginpar{$\K_u$}
\begin{equation}
  \label{eq:def:WQSym*:Ku}
\K_u := \N_u * E_n^{[r]} \quad (r=\ell(\ev(u))=\max(u))
\end{equation}
which will be a simultaneous eigenbasis for the modified Adams operations
$\psi_k(F) := F*(\sigma_1)^k$.

The closed expressions given in $\Sym$ can be readily extended to $\WQSym^*$
thanks to the following lemma.

\begin{lemma}\label{lem:perm}
Define a right action of $\SG_n$ on $\WQSym_n^*$ by
\begin{equation}
\N_u\cdot\sigma := \N_{u\sigma},\ \text{where}\ u\sigma = u_{\sigma(1)}u_{\sigma(2)}\cdots u_{\sigma(n)}.
\end{equation}
Then, for any $I\vDash n$,
\begin{equation}
\N_{u\sigma}*S^I = (\N_u*S^I)\cdot \sigma.
\end{equation}
\end{lemma}
\Proof If $S^I$ contains $\N_v$, it contains the $\N_{v\tau}$ for all permutations $\tau$,
and
\begin{equation}\label{eq:permpack}
\pack{u\tau\choose v} = \pack{u\choose v\tau^{-1}}\cdot\tau.
\end{equation}
\qed

For example, with $u=111122$, $v=212211$, $\tau=451623$,
we have $u\tau = 121211$, $v\tau^{-1}=211212$, 
$\pack{121211\choose 212211} =232411$, $\pack{111122\choose 211212}=211234$,
and $211234\tau = 232411$.

This implies\footnote{
Incidentally, this also proves the existence of the descent algebra. If one denotes by $\phi$
the dual of the inclusion map $\FQSym\rightarrow\WQSym$, which is given by
$\phi(\N_u)=\F_{\std(u)}$, then $\phi(S^I*\N_v) = \phi(S^I)\cdot\std(v)$, so that
$\phi(S^I*\N_v)=\phi(S^I)*\phi(\N_v)$.
} 
that $(f\cdot\sigma)*g =(f*g)\cdot\sigma$ for all $f\in\WQSym^*_n,
g\in\Sym_n$ and $\sigma \in\SG_n$.

Hence,
\begin{equation}
\psi_k(\N_u) = \sum_{v\ge u}\beta_k(v,u) \N_v,
\end{equation}
where $\beta_k(v,u):=\beta_k(\ev(v),\ev(u))$, and
\begin{equation}
  \label{eq:Ku=summuNv}
\K_u =\sum_{v\ge u}\frac{(-1)^{\ell(\ev(v))-\ell(\ev(u))}}{\ell(\ev(v),\ev(u))}\N_v,
\end{equation}
where the order on packed words $u,v$ is the {\it reverse} refinement order on
the corresponding set compositions $\sigma,\pi$ from §\ref{ssec:OST:poset}.

Given a spreadability system and a sequence $(X_i)$ of random variables, and defining the linear map 
$\hat\varphi_{\mathcal{X}}$ as above by $\hat\varphi_{\mathcal{X}}(\N_u)=\varphi_u(X_1,X_2,\ldots,X_n)$, the cumulants are
given by
\begin{equation}
K_u(X_1,X_2,\ldots,X_n)=\hat\varphi_{\mathcal{X}}(\K_u)
\end{equation}
and we recover the relations between moments
and cumulants of \cite[Th. 4.7 and 4.8]{HL}.

By construction,
\begin{equation}
\K_u*\S_n(N\alphA) = N^{\ell(\ev(u))}\K_u
\end{equation}
which is \cite[Th. 4.14]{HL}.
Indeed, according to \cite[Definition 4.1]{HL},
$$
\varphi_\pi(N.X_1,N.X_2,\ldots,N.X_n)=\hat\varphi_{\mathcal{X}}(\N_u*S_n(N\alphA))
$$
where $\pi=\oker(u)$.

Adapting the argument given at the end of \cite{NCSF7}, one can prove the
following extension of the ``splitting formula''
  \eqref{eq:Sym:internal-coprod}:
\begin{proposition}
If $f_1,f_2,\ldots,f_r\in\WQSym^*$ and $g\in\Sym$, then
\begin{equation}
(f_1f_2\cdots f_r)*g = \mu_r[(f_1\otimes\cdots\otimes f_r) *_r\Delta^rg].
\end{equation}
\end{proposition}
The same argument as in  \eqref{eq:mult} proves then the following product rule for the $\K_u$:
\begin{equation}
\N_u\N_v = \sum_w c_{u,v}^w \N_w\ \Rightarrow \ \K_u\K_v =  \sum_w c_{u,v}^w \K_w.
\end{equation}
That is, $\N_u\mapsto \K_u$ is an algebra automorphism.

\section{Partial cumulants}
The defining formula  \eqref{eq:phipi=sumksig}
of the cumulants can be inverted using the M\"obius function
\eqref{eq:moebiusmu}, but neither this nor the formula of Cartier and Good
\eqref{eq:omegacom} are suitable
for the efficient calculation of cumulants of higher orders.
In the case of exchangeability systems recursive formulas are available
which are more adequate for this purpose; see \cite[Proposition~3.9]{Leh04}.
In the classical case, the recursion reads as follows:
\begin{equation}
K(X_1,X_2,\dots,X_n) = 
\IE X_1X_2\dotsm X_n - \sum_{\substack{A\subsetneqq [n]\\ 1\in A}} K_{\abs{A}}(X_i:i\in A) \IE\prod_{j\in A^c}X_j
\end{equation}
and in the univariate case it specifies to
\begin{equation}
\kappa_n = m_n - \sum_{k=1}^{n-1}\binom{n-1}{k-1}\kappa_km_{n-k} 
\end{equation}
which is just a rearrangement of the familiar formula
  \eqref{eq:classicalcumulantrecurrence}.

In the free case it specifies
to the free Schwinger-Dyson equation \cite{MingoSpeicher:2013:schwinger}
and
from the point of view of combinatorial Hopf algebras this
has been recently considered under the name of ``splitting process''
\cite{EbrahimiFardPatras:2016:splitting}.

Turning to our general setting we note that already in the case of monotone
probability we lack a simple recursive formula; 
however Hasebe and Saigo~\cite{HS11b,HS11a} found
a good replacement in terms of differential equations.
This was generalized to spreadability systems in terms of
partial cumulants introduced in \cite[§6]{HL}
which are the images by $\hat\varphi_{\mathcal{X}}$ of some interesting elements of
$\WQSym^*$ which we shall study in this section.

We start with the analogous questions in $\Sym$, which are simpler and imply easily
the general results in $\WQSym^*$.
 

Let $T=(t_1,\ldots,t_r)$ be a sequence of binomial elements (scalars), so that the noncommutative
symmetric functions of $t_j\alphA$ are defined by
\begin{equation}
\sigma_1(t_j\alphA) := \sigma_1(\alphA)^{t_j}
\end{equation}
and the analogs of the formal multivariate moments \cite[Eq. (6.2)]{HL} are
\begin{equation}
S^I(T;\alphA) := S_{i_1}(t_1\alphA)\cdots S_{i_r}(t_r\alphA).
\end{equation}
The (generic) cumulants are thus
\begin{equation}
K_I = \left.\frac{\partial^r}{\partial t_1   \partial t_2\cdots \partial t_r}\right|_{T=(0,\ldots,0)}S^I(T;\alphA)
\end{equation}
Imitating \cite[Def. 6.1]{HL}, we define the partial cumulants  as
\begin{equation}
K_{I;j}^{(t_1,\ldots,t_{j-1},1,t_{j+1},\ldots,t_r)} := \left.\frac{\partial}{\partial t_j}\right|_{t_j=0}S^I(T;\alphA).
\end{equation}
Recall that $\sigma_1 = \exp(\phi)$, where $\phi=\sum_{n\ge 1}{\frac1n\Phi_n} = E^{[1]}$. Therefore,
\begin{equation}
\left.\frac{\partial}{\partial t}\right|_{t=0}\sigma_1(\alphA)^t = \phi,
\end{equation}
so that
\begin{equation}
\left.\frac{\partial}{\partial t_j}\right|_{t_j=0}S^I(T;\alphA)
 = S_{i_1}(t_1\alphA)\cdots S_{i_{j-1}}(t_{j_1}\alphA)\phi_{i_j} S_{i_{j+1}}(t_{j+1}\alphA)\cdots S_{i_r}(t_r\alphA),
\end{equation} 
and inserting
\begin{equation}
\phi_n = \sum_{H\vDash n}\frac{(-1)^{\ell(H)-1}}{\ell(H)}S^H =:  \sum_{H\vDash n}\tilde\mu(H,(n))S^H ,
\end{equation}
we obtain the analog of \cite[Prop. 6.2]{HL}
\begin{equation}
K_{I;j}^{(t_1,\ldots,t_{j-1},1,t_{j+1},\ldots,t_r)} 
= \sum_{K\vDash i_j}S^{I'HI''}((T',1^{\ell(H)},T'');\alphA)\tilde\mu(H,(n)), 
\end{equation}
where $T'=(t_1,\ldots,t_{j-1})$, $T''=(t_{j+1},\ldots,t_r)$ and $1^p$ means $(1,\ldots,1)$ ($p$ times).

Now,
\begin{equation}
\frac{\partial}{\partial t_j}S^I(T;\alphA)
\end{equation} 
is the homogeneous part of degree $n$ in
\begin{align}
S_{i_1}(t_1\alphA)\cdots S_{i_{j-1}}(t_{j-1}\alphA)  
\frac{\partial}{\partial t_j}e^{t_j\phi}
S_{i_{j+1}}(t_{j+1}\alphA)\cdots S_{i_r}(t_r\alphA)\\
=S_{i_1}(t_1\alphA)\cdots S_{i_{j-1}}(t_{j-1}\alphA)  
\phi\sigma_1^{t_j}
S_{i_{j+1}}(t_{j+1}\alphA)\cdots S_{i_r}(t_r\alphA)\\
=
S_{i_1}(t_1\alphA)\cdots S_{i_{j-1}}(t_{j-1}\alphA)  
\sigma_1^{t_j}\phi
S_{i_{j+1}}(t_{j+1}\alphA)\cdots S_{i_r}(t_r\alphA),
\end{align} 
the last two expressions being respectively equal to
\begin{equation}
\sum_{a=1}^{i_j} K_{(i_1,\ldots,i_{j-1},a,i_j-a,i_{j+1},\ldots,i_r);j}^{(t_1,\ldots,t_{j-1},1,t_j,\ldots,t_r)}
\end{equation}
and
\begin{equation}
\sum_{a=1}^{i_j} K_{(i_1,\ldots,i_{j-1},i_j-a,a,i_{j+1},\ldots,i_r);j+1}^{(t_1,\ldots,t_{j-1},t_j,1,\ldots,t_r)}
\end{equation}

These relations are then extended to $\WQSym^*$ by defining $\N_u(T)$ in such a way that
\begin{equation}
S^I(T;\alphA) = \sum_{\ev(u)=I}\N_u(T),
\end{equation}
which implies that (cf. \cite[Th. 4.5]{HL})
\begin{equation}
\N_u(T) = \sum_{v\ge u}\beta_T(v,u)\N_v.
\end{equation}

\begin{example}{\rm With $I=(2,2,1)$ we have
\begin{align}
S^{221}(T;\alphA)&= \left(t_1S_2+{t_1\choose 2}S^{11}\right)\left(t_2S_2+{t_2\choose 2}S^{11}\right)t_3S_1\\
&= t_1t_2t_3S^{221}+t_1{t_2\choose 2}t_3S^{2111}+{t_1\choose 2}t_2t_3S^{1121}+{t_1\choose 2}{t_2\choose 2}t_3S^{11111},
\end{align}
so that
\begin{align}
K_{(221);2}^{(t_1,1,t_3)}& =  \left(t_1S_2+{t_1\choose 2}S^{11}\right)(S_2-\frac12 S^{11})t_3S_1\\
&= t_1t_3S^{221}-\frac12t_1t_3s^{2111}+{t_1\choose 2}t_3S^{1121} -\frac12{t_1\choose 2}t_3S^{11111}. 
\end{align}
Similarly,
\begin{equation}
S^{2111}(T;\alphA)= \left(t_1S_2+{t_1\choose 2}S^{11}\right)t_2t_3t_4S^{111}
\end{equation}
and
\begin{equation}
K_{(2111);3}^{(t_1,t_2,1,t_4)}= \left(t_1S_2+{t_1\choose 2}S^{11}\right)t_2t_4S^{111}.
\end{equation}
We have indeed
\begin{align}
K_{(221);2}^{(t_1,1,t_3)}+K_{(2111);3}^{(t_1,t_2,1,t_3)} 
&= t_1t_3S^{221}+t_1t_3(t_2-\frac12)S^{2111}+{t_1\choose 2}t_3S^{1121}+{t_1\choose 2}(t_2-\frac12)t_3S^{11111}\nonumber\\
&=\frac{\partial}{\partial t_2}S^{221}(T;\alphA).
\end{align}
In $\WQSym^*$, this would translate for the set composition $\pi=24|15|3\leftrightarrow u=21312$ as
\begin{equation}
\N_{21312}(T;\alphA) =
 t_1t_2t_3\N_{21312}+t_1{t_2\choose 2}t_3\N_{21413}
+{t_1\choose 2}t_2t_3\N_{31423}+{t_1\choose 2}{t_2\choose 2}t_3\N_{31524},
\end{equation}
and
\begin{equation}
\N_{21413}(T;\alphA) =
 t_1t_2t_3t_4\N_{21413}+{t_1\choose 2}t_2t_3t_4\N_{31524},
\end{equation}
and the partial cumulants $K_{24|15|3;15}^{(t_1,1,t_3)}$, $K_{24|1|5|3;5}^{(t_1,t_2,1,t_4)}$
would be respectively
\begin{equation}
K_{21312;2}^{(t_1,1,t_3)} = t_1t_3\N_{21312}-\frac12t_1t_3\N_{21413}+{t_1\choose 2}t_3\N_{31423}-\frac12{t_1\choose 2}t_3\N_{31524}
\end{equation}
\begin{equation}
K_{21413;3}^{(t_1,t_2,1,t_4)} = t_1t_2t_4\N_{21413}+{t_1\choose 2}t_2t_4\N_{31524}
\end{equation}

}
\end{example}

\section{Left and right products with the Eulerian idempotents}

We have seen that the cumulant basis is given by the internal products
$\N_u*E^{[r]}$ where $r=\ell(\ev(u))=\max(u)$.

The aim of this section is to compute the internal products
\begin{equation}
\N_u*E^{[k]}\quad\text{and}\quad  E^{[k]}*\N_u
\end{equation}
for arbitrary $u$ and $k$.

Let $u$ be a packed word. Recall that $v$ is said to refine $u$ if for all
 $i<j$, $v_i>v_j \Longleftrightarrow u_i\ge u_j$ and $v_i=v_j\Longrightarrow u_i=u_j$.
In this case, we write $v\ge u$ or  $v\in \raff(u)$.
This is the usual notion of refinement on set compositions from
§\ref{ssec:OST:poset}:
each block of $u$ is a
union of \emph{consecutive} blocks of $v$. 

We shall say that $v$ is a \emph{weak refinement} of $u$, and write
$v\succeq u$ or $v\in \raffbis(u)$, if for all  $i<j$, $v_i=v_j\Longrightarrow u_i=u_j$.
On set compositions, this means that each block of $u$ is a union of blocks of $v$.

For example,
\begin{equation}
\begin{split}
& \raff(122) = \{122, 123, 132\}, \\
&\raffbis(122) = \{122,211,123,132,213,231,312,321\}.
\end{split}
\end{equation}

Thus, the word $w=\pack{u\choose v}$ is finer than $u$ and is such that
$w_i=w_j \Longrightarrow v_i=v_j$.

We shall compute the generating functions
\begin{align}
U_v(t) &= \sigma_1^t*\N_v = \sum_r t^r E^{[r]}*\N_v =  \sum_u {t\choose\max(u)}\N_u*\N_v\\
&=\sum_{w\in\raffbis(v)}\left( \sum_{u\in U(v,w)}{t\choose\max(u)}\right)\N_w,
\end{align}
where 
\begin{equation}
U(v,w)=\left\{u|\pack{u\choose v}=w\right\}
\end{equation}
and
\begin{align}
V_u(t) &= \N_u*\sigma_1^t = \sum_r t^r \N_u*E^{[r]} =  \sum_v {t\choose\max(v)}\N_u*\N_v\\
&=\sum_{w\ge u}\left( \sum_{v\in V(u,w)}{t\choose\max(v)}\right)\N_w,
\end{align}
where 
\begin{equation}
V(u,w)=\left\{v|\pack{u\choose v}=w\right\}.
\end{equation}


\subsection{A closed formula for $U_v(t)$}

The  pairs $(v,w)$ such that $U(v,w)$  is nonempty are those such that $w\in\raffbis(v)$.

\begin{proposition}
Define a set composition $u_0=u_0(v,w)$ as the one obtained by merging two consecutive blocks
$p',p''$ of $w$ if the blocks of $v$ containing $p'$ are strictly to the left of those containing $p''$, 
and iterating the process until no further blocks can be merged.

Then, $U(v,w)$ is the interval
\begin{equation}
U(v,w) = [u_0,w]\quad \text{(compositions finer than $u_0$ and coarser than $w$).}
\end{equation}
\end{proposition}

For example, let us compute
$U(13211, 15342)$. The corresponding set compositions are
$145|3|2$ and $1|5|3|4|2$. To construct $u_0$, one can
merge the second and third blocks of $w$, and also the fourth and the fifth ones,
which yields $u_0=1|35|24$.
Therefore,
\begin{equation}
U(13211, 15342) = [13232,15342] = \{13232, 14232, 14342, 15342\}.
\end{equation}

To compute $U(13211, 24315)$, the corresponding set compositions are
$145|3|2$ and $4|1|3|2|5$, we obtain $u_0=4|123|5$, so that
\begin{equation}
U(13211, 24315) = [22213,24315] = \{22213, 23214, 23314, 24315\}.
\end{equation}

Finally, with $U(13211, 13214)$, the set compositions are 
 $145|3|2$ and $14|3|2|5$,  $u_0=1234|5$,
and
\begin{equation}
U(13211, 13214) = [11112,13214] = \{11112, 12113, 12213, 13214\}.
\end{equation}

\Proof It follows from the definition of $\pack{u\choose v}$ that if two elements
$i,j$ are in the same block of $w$, then, they must also be in the same block in $u$ and in $v$, 
and otherwise $i$ is in a block strictly left to the block of $j$ in $v$ if and only if
either $i$ is in a block left to the block of $j$ in $u$, or $i$ and $j$ are in the same
block of $u$, and the block of $i$ if left to the block of $j$ in $v$.  

Thus, any $u\in U(v,w)$ must be obtained by merging some consecutive parts in $w$,
so that $u\le w$. Moreover, to get $w$ as $\pack{u\choose v}$, one can only merge
blocks satisfying the constraints mentioned in the definition of $u_0$, so that $u\ge u_0$.  
\qed

So, $U(v,w)$ is the interval $[u_0,w]$. This interval is clearly a Boolean lattice.
Moreover, if we set $a_0=\max(u_0)$ and $a=\max(w)$, the number of elements of
this lattice such that $\max(u)=k$ is ${a-a_0\choose a-k}$, so that
\begin{equation}
\sum_{u\in U(v,w)}
{t\choose\max(u)}
=
\sum_{k=a_0}^a {a-a_0\choose a-k} {t\choose k}
= {t+a-a_0\choose a},
\end{equation}
and finally
\begin{equation}\label{eq:U}
U_v(t) = \sum_w {t+\max(w)-a_0(v,w)\choose\max(w)}\N_w
\end{equation}
where $a_0(v,w) = \max(u_0(v,w))$.

In particular, the coefficient of $\N_w$ in $E^{[1]}*\N_v$ is the coefficient of
$t\N_w$ in $U_v(t)$, which is
\begin{equation} 
(-1)^{a_0-1}\frac{(a-a_0)!(a_0-1)!}{a!}= (-1)^{a_0-1}{\rm B}(a-a_0+1,a_0)
\end{equation}
where ${\rm B}$ is the Beta function. This is formula (7.3) of \cite{HL} for the Weisner coefficients.

\subsection{A closed formula for $V_u(t)$}

Let us now describe $V(u,w)$. 
Since the packing process commutes with the right action of the symmetric group (see \eqref{eq:permpack}),
we can apply to $u$
the smallest permutation $\sigma$  
such that $u\sigma$ is nondecreasing (i.e., $\sigma=\std(u)^{-1}$), so that
$\pack{u\sigma\choose v\sigma}=w\sigma$. We can therefore assume that $u$ is nondecreasing.

\begin{proposition}
Let $w^{(i)}=\pack(w_{j_1}w_{j_2}\cdots w_{j_p})$, where $\{j_1,\ldots,j_p\}=\{j|u_j=i\}$.
Then,
\begin{equation}
\sum_{v\in V(u,w)}\M_v = \M_{w^{(1)}}\M_{w^{(2)}}\cdots \M_{w^{(\max(u))}}.
\end{equation}
\end{proposition}

\Proof 	Note first that no relation is imposed between the letters of $v$ corresponding
to different letters of $u$. The only order  constraints are among places where $u$ has identical letters,
and these are the same as in the corresponding letters of $w$. This is precisely the definition of the convolution
on packed words, describing the product of the $\M$ basis. \qed

Since the map $\M_u\mapsto {t\choose u}$ is a character of $\WQSym$,
\begin{equation}
\sum_{v\in V(u,w)} \binom{t}{\max(v)}
= \prod_{i=1}^{\max(u)} \binom{t}{\max(w^{(i)})}.
\end{equation}
We have therefore
\begin{equation}
V_u(t) = \sum_{w\ge u}\prod_{i=1}^{\max(u)} \binom{t}{\max(w^{(i)})}\N_w.
\end{equation}

For example, take $u=1122$ and $w=2133$.
We have $w^{(1)}=21$ and $w^{(2)}=\pack(33)=11$.
The product is
\begin{equation}
\M_{21} \M_{11} = \M_{2111} + \M_{2122} + \M_{2133} + \M_{3122} + \M_{3211},
\end{equation}
and the set of  $v$ is
\begin{equation}
V(1122,2133) = \{2111, 2122, 2133, 3122, 3211\}.
\end{equation}
We can then see that
\begin{equation}
2\binom{t}{2} + 3\binom{t}{3} = \frac{t^2(t-1)}{2} = \binom t2 \binom t1,
\end{equation}
as claimed.

\subsection{Mixed cumulants}

The mixed cumulants described in §7 of \cite{HL} correspond to the elements $\K_u*\N_v=\N_u*E^{[\max(u)]}*\N_v$
of $\WQSym^*$.

We shall compute the generating series
\begin{equation}
\N_u * \sigma_1^t * \N_v=\sum_r \N_u*E^{[r]}*\N_v.
\end{equation}
This amounts to computing
\begin{equation}
\label{all}
 \sum_{u'\ge u}
           \prod_{i=1}^{\max(u)} \binom{t}{\max(u'^{(i)})} \N_{u'}*\N_v.
\end{equation}

Let us compute the coefficient of $\N_w$ in~\eqref{all}.
First, the packed words having a nonzero coefficient in this expansion
are those finer than $u$, which are weakly finer than $v$.
Let $w$ be such a word. The coefficient of $\N_w$ is obtained by summing the coefficients of all $u'\ge u$ such that
$\pack{u'\choose v}=w$. The set of those $u'$ is therefore the set
of packed words which are finer than $u$ and than 
$u_0(v,w)$, and that are coarser than $w$.

This set is an intersection of Boolean lattices, hence also a Boolean lattice.
Moreover, the subwords of $u$ consisting of identical letters cannot interfere
with each other, so that this lattice is in fact the product of the lattices
obtained by restricting $v$ and $w$ to positions where $u$ has identical letters.

Thus, the result is obtained by applying~\eqref{eq:U} 
on each set of positions where $u$ has identical letters since these
pieces are independent and $\N_{1^k}=S_n$ is neutral for $*$. \qed

The result is therefore
\begin{equation}
\N_u * \sigma_1^t * \N_v
 = \sum_{w\in W(u,v)}
     \prod_{i=1}^{\max(u)} \binom{t+\max(w^{(i)})-a(v^{(i)},w^{(i)})}{\max(w^{(i)})} \N_{w},
\end{equation}
where $a$ is as previously defined, and
$W(u,v)$ 
is the set of packed words finer than $u$ which have equal letters only at
positions where $u$ has equal letters, and
$v^{(i)}$ and $w^{(i)}$ 
are the subwords of $v$ and $w$ corresponding to the positions of the letter $i$ in $u$.

\medskip
\begin{example}{\rm Let $u=11122$ and $v=12234$.

Take $w=23145$. We form the packed words of the restrictions of    $v$ and
$w$ to the positions where $u$ has equal letters, which gives
for the first block $v^{(1)}=122$ and $w^{(1)}=231$, 
with a contribution of  ${t+1\choose 3}$  for this factor
(since $u_0^{(1)}=221$) and for the second,
$v^{(2)}=12$ and $w^{(2)}=12$, hence a contribution of ${t+1\choose 2}$ (since
$u_0^{(2)}=11$).

For $w=21354$, we have respective contributions of ${t+1\choose 3}$ and
${t\choose 2}$ since $u_0^{(1)}=212$ and $u_0^{(2)}=21$.

For $w=31245$,  we have respective contributions  of ${t\choose 3}$ and
${t+1\choose 2}$ since $u_0^{(1)}=312$ and $u_0^{(2)}=11$.

For  $w=31254$,  we have respective contributions  of ${t\choose 3}$ and
$\binom{t}{2}$ since $u_0^{(1)}=312$ and $u_0^{(2)}=21$.

Finally, for  $w=12243$,  we have respective contributions of ${t+1\choose 2}$ and
${t\choose 2}$ since $u_0^{(1)}=111$ and $u_0^{(2)}=21$.
}
\end{example}

\section{Miscellaneous remarks}

We have already noticed (see \eqref{eq:permpack}) that
\begin{equation}
\pack\gf {u\cdot \sigma}{v\cdot \sigma} = \pack\gf uv\cdot\sigma
\end{equation}
for any permutation $\sigma$.

Taking the smallest permutations sorting $u$ ($\sigma=\std(u)^{-1}$),
we can restrict to the case where $u$ is nondecreasing. With a second
permutation, we can restrict to the case where $v$ is nondecreasing
on positions where $u$ has equal letters.
Moreover, if
\begin{equation}
\pack\gf u{v_1} = \pack\gf u{v_2},
\end{equation}
then
\begin{equation}
\pack\gf {u'}{v_1} = \pack\gf {u'}{v_2}
\end{equation}
for all $u'$ finer than $u$. In the latter argument, it is thus possible
to restrict to the case where $u$ is a shifted concatenation of packed words,
one per block of equal letters of $u$.

\section{A generalization of Goldberg's formula}

The product of exponentials
\begin{equation}
g =  e^{a_1}e^{a_2}\cdots = \sum_{u\ {\rm packed}}\frac1{u!}\M_u(\alphA), \quad (u! := \prod_i |u|_i!)
\end{equation}
is naturally an element of the completion $\widehat{\WQSym}$ of $\WQSym$. It is grouplike for the coproduct of $\WQSym$, which in this case
coincides with the coproduct of $\KK\llangle \alphA\rrangle$. 
The Hausdorff series
\begin{equation}
H(a_1,a_2,\ldots) = \log g = \sum_u c_u \M_u
\end{equation}
is thus also an element of $\widehat{\WQSym}$. The coefficient $c_u$ can be computed by Goldberg's formula, of which
new proofs respectively based on arguments from combinatorial Hopf algebras and from noncommutative probability
have been recently given in \cite{FPT} and \cite{HL}.  

The homogeneous component $H_n$ of $H$ is the image of $g_n$ by the first Eulerian idempotent $e_n^{[1]}$, {\it aka} Solomon's
idempotent. More generally, the image of $g_n$ by $e_n^{[k]}$ is the coefficient of $t^k$ in $g^t$. Define polynomials $c_u(t)$ by
\begin{equation}
g^t =  (e^{a_1}e^{a_2}\cdots)^t =:\sum_u c_u(t)\M_u(\alphA).  
\end{equation}
Since $g$ is grouplike, we can, as in \cite{FPT}, define an injective morphism of Hopf algebras 
\begin{align}
\varphi:\ & \Sym\longrightarrow \WQSym\\
&           S_n\longmapsto g_n
\end{align}
Then, $g^t =\varphi(\sigma_1^t)$ and
\begin{align}
c_u(t) &= \< \N_u,\varphi(\sigma_1^t)\> =\<\varphi^\dagger(\N_u),\sigma_1^t\>\\
&= \sum_I {t\choose \ell(I)}\<\varphi^\dagger(\N_u),S^I\>
\end{align}
where $\varphi^\dagger:\ \WQSym^*\rightarrow QSym$ is the adjoint map, and $\<\cdot,\cdot\>$
is the duality bracket.

Now, recall that the coproduct of $\N_u$ is given by \cite{NT06}
\begin{equation}
\Delta\N_u = \sum_{u=u_1\cdot u_2}\N_{\pack(u_1)}\otimes \N_{\pack(u_2)}
\end{equation}
(deconcatenation). We can omit the packing operation in this formula
if we make the convention that $\N_w=\N_u$ if $u=\pack(w)$. Then,
since $\varphi$, and hence also $\varphi^\dagger$ are morphisms of
Hopf algebras, for a composition $I=(i_1,\ldots,i_r)$,
\begin{equation}
\<\varphi^\dagger(\N_u),S^I\>=\prod_{k=1}^r\<\varphi^\dagger(\N_{u_k}),S^{i_k}\>
\end{equation}
where $u=u_1u_2\cdots u_r$ with $|u_k|=i_k$ for all $k$. Moreover, this is nonzero
only if all the $u_k$ are nondecreasing, in which case the result is
$1/(u_1!\cdots u_r!)$.

Let $u=u_1\cdots u_r$ be the minimal factorization of $u$ into nondecreasing words
({\it i.e.}, such that the last letter of each $u_i$ is strictly greater than the first one
of $u_{i+1}$), and let $I=(|u_1|,\ldots,|u_r|)$.

Let also $J=(j_1,\ldots,j_s)$ be the composition obtained by factoring $u$ into maximal blocks of 
identical letters,
\begin{equation}
u = b_1^{j_1}b_2^{j_2}\cdots b_s^{j_s}.
\end{equation}
Then, it is worth noticing that $c_u(t)$ is the $(t,\E)$ specialization of
\begin{equation}
c_u(X;\alphA) = \sum_{K\ge I} M_K(X)S^{K\vee J} = \sum_{H\ge J}\left(\sum_{K\vee J=H}M_K(X)\right)S^H(\alphA) \quad\in QSym\otimes \Sym,
\end{equation}
where $\vee$ denotes the join in the lattice of compositions of $n$.
For a binomial element $t$, $M_K(t)={t\choose\ell(K)}$, and the (virtual) exponential alphabet $\E$ is defined
by $S_k(\E)=\frac1{k!}$.
The set $\{K|K\vee J=H\}$ is a Boolean lattice: it is the interval $[I\vee L,H]$ where  $\Des(L)=\Des(H)\backslash\Des(J)$.
Therefore,
\begin{equation}\label{eq:cu}
c_u(t) = \sum_{H\ge J}{t+\ell(J)-\ell(I)\choose\ell(H)}S^H(\E).
\end{equation}

For example, with $u=113223$, we have $I=33$, $J=2121$,
\begin{align}
c_u(X,\alphA)=& (M_{33}+M_{213}+ M_{321}  + M_{2121})S^{2121}\nonumber\\
& + (M_{123} + M_{1113} + M_{1221} + M_{11121})S^{11121}\nonumber\\
& + (M_{312} + M_{3111} + M_{2112} + M_{21111})S^{21111}\nonumber\\
& +( M_{1212} + M_{12111} + M_{11112} + M_{111111})S^{111111} 
\end{align}
whose  $(t,\E)$-specialization  reduces, by binomial convolution,  to
\begin{equation}
{t+2 \choose 4}\frac1{2!1!2!1!} +{t+2 \choose 5 }\frac1{1!1!1!2!1!} +{t+2 \choose 5}\frac1{2!1!1!1!1!} +{t+2 \choose 6 }\frac1{1!1!1!1!1!1!} 
\end{equation}
This is the image of the polynomial
\begin{equation}
\left[\frac{t(2t^2+t)}{2}\right]^2 =: E^{2121}(t)
\end{equation}
under the linear substitution $t^k\mapsto {t+2\choose k}$.

Define a linear map ${\mathcal F}_s:\ t^k\mapsto {t+s\choose k}$, and associate to a composition $J$
the product of normalized Eulerian polynomials
\begin{equation}
E^J(t) := \prod_{k=1}^{\ell(J)}\frac{t E_{j_k}(t,t+1)}{j_k!}
\end{equation} 
as in Goldberg's formula \cite[Theorem 3.11]{Reu},
where
\begin{equation}
E_n(x,y)=\sum_{\sigma\in\SG_n}x^{d(\sigma)}y^{n-d(\sigma)},
\end{equation}
and $d(\sigma)$ is the number of descents of $\sigma$.  

\begin{theorem}
The $t$-Goldberg coefficient $c_u(t)$ is
\begin{equation}\label{eq:gg}
c_u(t) = {\mathcal F}_{s-r}(E^J(t)),
\end{equation}
where $r=\ell(I)$, $s=\ell(J)$, and $I,J$ are the compositions recording  respectively the lengths of the nondecreasing runs
and of the maximal blocks of identical letters of $u$.
\end{theorem}

The coefficient of $t$ in ${\mathcal F}_{s-r}(t^k)$ is
\begin{equation}
\begin{split}
(-1)^{k-s+r-1}\frac{(s-r)!(k-s+r-1)!}{k!}&=(-1))^{k-s+r-1}{\rm B}(s-r+1,k-s+r)\\
&=\int_{-1}^0t^k\cdot t^{r-1}(1+t)^{s-r}\frac{dt}{t^s}
\end{split}
\end{equation}
(since $k\ge s$ and $r\ge 1$), which applied to \eqref{eq:gg} yields the classical form of Goldberg's formula:
\begin{equation}
c_u = [t]c_u(t) = \int_{-1}^0 t^{r-1}(1+t)^{s-r} \prod_{k=1}^{s}\frac{ E_{j_k}(t,t+1)}{j_k!}dt.
\end{equation}
Note that $r-1$ is the number of descents of $u$, and $s-r$ its number of rises.

\section{Mixed cumulants and Goldberg coefficients}

It remains to explain the occurrence of Goldberg coefficients in the expansion of mixed
cumulants on the cumulant basis.

As before, we shall work with the $t$-analogues.

Let $u$ be a packed word.
Define as above two compositions $I(u)$ and $J(u)$ recording the lengths
of the maximal nondecreasing factors and of the maximal blocks of identical letters of $u$.
For example,
\begin{equation}
I(31121) = (1,3,1) \text{\quad and \quad} J(31121)=(1,2,1,1).
\end{equation}
Note that $J(u)\ge I(u)$.

Following \cite{HL}, with a pair of words such that $v\succeq u$, we associate their
\emph{refinement word} $m(u,v)$ defined as follows: $|m(u,v)|=\max(v)$, and 
$m_i=u_j$ whenever $v_j=i$. For example,
\begin{equation}
m(12113,41223)= 2131,\
m(111123,353241) = 31121.
\end{equation}

We can now reformulate \eqref{eq:U} as
\begin{equation}\label{eq:newU}
U_v(t) = \sigma_1^t * \N_v
 = \sum_{ \ogf{w\in\raffbis(v)}{m=m(v,w)} } \binom{t+\ell(J(m))-\ell(I(m))}{|J(m)|} \N_w.
\end{equation}
For example,
\begin{equation}
\begin{split}
U_{112}(t) &=
\binom{t+1}{2} \N_{112} + \binom{t}{2} \N_{221} +
\binom{t+1}{3} \N_{123} + \binom{t+1}{3} \N_{132} \\
&+
\binom{t+1}{3} \N_{213} + \binom{t}{3} \N_{231} +
\binom{t+1}{3} \N_{312} + \binom{t}{3} \N_{321}.
\end{split}
\end{equation}
Computing $U_{1122}(t)$, one finds a sum of 38 terms, corresponding
to the packed words of which the first two or the last two letters can be identical.
There are 7 different coefficients in this expansion, given below together
with their associated words:
\begin{equation}
\begin{split}
\binom{t+2}{4} &: \{1324, 1342, 3124, 3142\} \\
\binom{t+1}{4} &: \{1234, 1243, 1423, 1432, 2134, 2143, 2314, 2341, 2413, \\
              &\ \ \ 2431, 3214, 3241, 4123, 4132, 4213, 4231\} \\
\binom{t+1}{3} &: [1123, 1132, 1233, 1322, 2133, 2213, 2231, 3122] \\
\binom{t+1}{2} &: \{1122\} \\
\binom{t}{4}   &: \{3412, 3421, 4312, 4321\} \\
\binom{t}{3}   &: \{2311, 3211, 3312, 3321\} \\
\binom{t}{2}   &: \{2211\}
\end{split}
\end{equation}

Let us now establish the equivalence of \eqref{eq:U} and \eqref{eq:newU}. 
We have to prove that for all
 $v$ and $w$
such that $w\in\raffbis(v)$ :
\begin{equation}
\binom{t+\max(w)-a_0(v,w)}{\max(w)} =
\binom{t+\ell(J(m))-\ell(I(m))}{|J(m)|}.
\end{equation}
First of all, it is clear, by definition, that $|J(m)|=\max(w)$.
Moreover, $\max(w)-a_0(v,w)$ is the number of blocks of $w$ (regarded as a
set composition) whose all elements are all strictly to the left in $v$ of all
elements of the next block. This amounts precisely to merging two parts of $J(m)$
whenever the corresponding values in $m$ values are in increasing order. 
Thus, $\max(w)-a_0(v,w)=\ell(J(m))-\ell(I(m))$.

Consider now the expansion of $U_v(t)$ on the basis $\K_w$.
Recall the transition matrices between $\N$ and $\K$. 

\begin{equation}
\K_v = \sum_{w\in\raff(v)} 
        \frac{(-1)^{\max(w)-\max(v)}}{\pev(m(u,v))} \N_w,
\end{equation}
\def\pev{{\pi_{\rm ev}}}
where $\pev(m)=\prod_i|m|_i$.

Notice that by definition
\begin{equation}
\K_v = [t^{\max(v)}] \N_v *\sigma_1^t.
\end{equation}
Indeed, to obtain the correct power of $t$, one must select in each
$\binom{t}{M}$ the coefficient of  $t$, which is $\frac{(-1)^{M-1}}{M}$.

In the other direction,
\begin{equation}
\N_v = \sum_{w\in\raff(v)}
        \frac{1}{\pev(m(u,v))!} \K_w,
\end{equation}
where $\pev(m)!=\prod_i|m|_i!$.

Ths yields for the expansion of $U_v(t)$ on the basis $\K$
\begin{equation}
U_v(t) = \sum_{x\in\raffbis(v)} c_{m(v,x)}(t) \K_x,
\end{equation}
where $c_{m(v,x)(t)}$ is defined in Eq.~\eqref{eq:cu}.

Indeed, to go from the expression of $U_v$ on the $\N_w$ to its
expression on the $\K_x$, notice that 
the $\N_w$ contributing to a given $\K_x$ are those such that $w\le x$ and
$w\succeq v$.
This is a Boolean lattice, an interval $[x_0,x]$ where $x_0$ is obtained by
merging consecutive blocks of $x$ whenever these blocks are contained in a
block of $v$.

We must therefore evaluate
\begin{equation}
\sum_{w\in [x_0,x]} \frac{1}{\pev(m(w,x))!}
\binom{t+\ell(J(m(v,w)))-\ell(I(m(v,w)))}{|J(m(v,w))|}.
\end{equation}
For a $w\in[x_0,x]$, it is easy to compute $m(v,w)$. The set of those
$m(v,w)$ is the set of words obtained from $m(v,x)$ by replacing blocks
$i^r$ of consecutive identical letters by smaller blocks $i^p$, $1\le p\le r$.
On this set,
$\ell(J(m))-\ell(I(m))$ is constant, 
 $H=\ev(m(w,x))$
runs over the set of compositions finer than $\ev(m(x_0,x))$,
$1/\ev(m(w,x))!=S^H(\E)$ and $|J(m(v,w))|=\max(w)=\ell(H)$.

For example, let $v=111123$ and $x=356241$.
Then, $x_0=244231$.
The contribution to $\K_x$ in $U_v(t)$ 
come from the four words
\begin{equation}
\{ 244231, 245231, 355241, 356241\}.
\end{equation}
Their words $m(v,w)$ are respectively $3121$, $31211$, $31121$ and $311211$.
The  $m(w,x)$ are $122344$, $122345$, $123455$ and $123456$.
The sum of all contributions is
\begin{equation}\label{ex:gm}
\frac{1}{4}\binom{t+1}{4} + \frac{1}{2}\binom{t+1}{5} +
\frac{1}{2}\binom{t+1}{5} + \binom{t+1}{6}.
\end{equation}
This is the same as $c_u=c_{311211}(t)$.
Indeed,  $I(u)=(1,3,2)$ and $J(u)=(1,2,1,2)$, so that we get
\begin{equation}
\sum_{H\geq 1212} \binom{t+1}{\ell(H)} S^H(\E)
\end{equation}
which yields the same expression as \eqref{ex:gm}.

Finally, one can give the general formula for
$\N_u * \sigma_1^t * \N_v$ on the $\K_w$.
This is immediate, since the transition from $\N$ to $\K$
does not change the structure of the product. We get
\begin{equation}
\N_u * \sigma_1^t * \N_v
 = \sum_{w\in W(u,v)}
     \prod_{i=1}^{\max(u)} c_{m^{(i)}}(t) \K_w,
\end{equation}
where $c_u(t)$ are the  $t$-Goldberg coefficients, and
$m^{(i)}=m(v^{(i)},w^{(i)})$.


\footnotesize
\bibliographystyle{amsalpha}
\bibliography{spreading}

\end{document}